\RequirePackage{fix-cm}
\documentclass[smallextended]{svjour3}       % onecolumn (second format)
\smartqed  % flush right qed marks, e.g. at end of proof
\usepackage{graphicx}
\usepackage{amsmath}
\usepackage{amssymb}
\usepackage{theorem}
\usepackage{xspace}
\usepackage{epsfig}
\usepackage{hhline}
\usepackage[utf8]{inputenc}    % pour les accents 
\def\hh{ \hspace*{0.5cm}}
\def\argmax{\mathop{\mathrm{arg\,max}}}

\def\hh{ \hspace*{0.5cm}}

\def\eb{\textrm{\mathversion{bold}$\mathbf{\beta}$\mathversion{normal}}}  
\def\ebo{\textrm{\mathversion{bold}$\mathbf{\beta^0}$\mathversion{normal}}}
\def\el{\textrm{\mathversion{bold}$\mathbf{\lambda}$\mathversion{normal}}}  
\def\ephi{\textrm{\mathversion{bold}$\mathbf{\phi}$\mathversion{normal}}}

\def\eR{I\!\!R}
\def\eE{I\!\!E}
\def\eP{I\!\!P}
\def\e1{1\!\!1}

\def\ef{\mathbf{\overset{.}{f}}}
\def\eff{\mathbf{\overset{..}{f}}}

\def\eg{\mathbf{{g}}}
\def\egg{\mathbf{\overset{.}{g}}}
\def\eV{\mathbf{{V}}}

\def\eVV{\mathbf{{V}}}
\def\ev{\mathbf{{v}}}
\def\eD{\mathbf{{D}}}
\def\eM{\mathbf{{M}}}

\def\ez{\mathbf{{z}}}
\def\XX{\textrm{\mathversion{bold}$\mathbf{X}$\mathversion{normal}}}
\def\xx{\textrm{\mathversion{bold}$\mathbf{x}$\mathversion{normal}}}
\def\ell{\textrm{\mathversion{bold}$\mathbf{\overset{.}{\lambda}}$\mathversion{normal}}}
\newcommand{\R}{\mathbb{R}}
\theoremstyle{plain}

\newcommand{\Var}{\mathbb{V}\mbox{ar}\,}
\def\argmax{\mathop{\mathrm{arg\,max}}}
 
\newtheorem{Remark}{Remark} %[section]
\newtheorem{Corollary}{Corollary} 
\begin{document}

\title{Empirical likelihood test in a posteriori  change-point nonlinear model}




\author{Gabriela Ciuperca         \and
        Zahraa Salloum 
        }


\institute{ G.Ciuperca - Z. Salloum \\
Universit\'e de Lyon, Universit\'e Lyon 1, CNRS, UMR 5208, Institut Camille Jordan, \at
             Bat.  Braconnier, 43, blvd du 11 novembre 1918, F - 69622 Villeurbanne Cedex, France \\
                         \email{Gabriela.Ciuperca@univ-lyon1.fr}           \\ \\
            % \emph{Present address:} of F. Author  %  if needed
           %\and
       %   Universit\'e de Lyon, Universit\'e Lyon 1, CNRS, UMR 5208, Institut Camille Jordan, \at
         %    Bat.  Braconnier, 43, blvd du 11 novembre 1918, F - 69622 Villeurbanne Cedex, France \\
       %    Tel.: 33(0)4.72.44.62.19\\
        %   Fax: 33(0)4.72.43.16.87\\
    Z. Salloum \\
            \email{salloum@math.univ-lyon1.fr}      
}

\maketitle

\begin{abstract}
\hh In this paper, in order to test whether changes have occurred  in  a nonlinear parametric regression, we propose a nonparametric method based on the empirical likelihood. Firstly, we test the null hypothesis of no-change against the alternative of one change in the regression parameters. Under null hypothesis, the consistency and the convergence rate of the regression parameter estimators are proved. The asymptotic distribution of the test statistic under the null hypothesis is obtained, which allows to find the asymptotic critical value. On the other hand, we prove that the proposed test statistic has the asymptotic power equal to 1. These theoretical results allows find a simple test statistic, very useful for applications.  The epidemic  model, a particular model  with two change-points under the alternative hypothesis, is also studied. Numerical studies by Monte-Carlo simulations show the performance of the proposed test statistic, compared to an existing method in literature.

\end{abstract}

\section{Introduction}
\label{intro}
We consider a classical model of parametric nonlinear regression : 
\begin{equation}
\label{eqq1}
Y_i=f(\XX_i; \eb)+\varepsilon_i, \qquad i=1, \cdots, n,
\end{equation}
where a possible change in the regression parameters could occurs. This is called, change-point problem. \\
Change-point detection problems fall in two categories. The first type is \textit{a posteriori}: after that  the $n$ all observations are realized, we study if,  a certain moment $k \in \{2, \cdots , n-1 \}$, the model (parameter $\eb$, to be more precise) is changed :
\begin{equation}
\label{eqq2}
Y_i=\left\{
\begin{array}{ccl}
f(\XX_i;\eb_1)+\varepsilon_i &  & i=1, \cdots, k \\
f(\XX_i;\eb_2)+\varepsilon_i &  & i=k+1, \cdots, n.
\end{array}
\right.
\end{equation}
The second type of change-points model is sequential (\textit{a priori}), where the change detection  is performed  in  real time. If in the first $k$ observations no change in the parameter regression has occurred, at observation $k+1$ we test that there is no change in the model:  $Y_i=f(X_i; \eb)+\varepsilon_i$, for all $i=1, \cdots , k+1$, against the hypothesis that the model has the form :
\begin{equation}
\label{eqq3}
\begin{array}{lcl}
Y_i=f(\XX_i;\eb)+\varepsilon_i &  & \textrm{for } i=1, \cdots, k \\
Y_{k+1}=f(\XX_{k+1};\eb^*)+\varepsilon_{k+1}, &  & 
\end{array}
\end{equation}
with $\eb \neq \eb^*$.  \\
\hh In this paper, we consider a posteriori  change-point problem. \\ 
\hh For the two types of problems, the number of publications in the last years is every extensive. Let us mention some  references concerning  the sequential change-point problem. If the function $f$ is linear, $f(\xx, \eb)=\xx^t \eb$, in the papers \cite{Horvath:Huskova:Kokoszka:Steinebach:04}, \cite{Huskova:Kirch:12},  the CUSUM method is used to find a test statistic for detecting the presence or absence of a change. The results have been generalized by \cite{Ciuperca:13a} for a nonlinear model. We can also mention the  papers \cite{Lai:Xing:10}, \cite{Mei:06}, \cite{Neumeyer:Van Keilegom:09}  for the sequential detection of a change-point. \\
\hh For a posteriori  change-point problem, in order to detect a change-point presence, model (\ref{eqq1}) is tested against model (\ref{eqq2}). The non-identifiability of model under the null hypothesis makes  classical test techniques unusable. In most articles in the literature, the authors propose criteria: see for example \cite{Nosek:10}, \cite{Ciuperca:11}, \cite{Wu:08}. Various hypothesis tests have been proposed only for the  linear models. The likelihood-ratio test method is used in \cite{Bai:99} and \cite{Lee:Seo:Shin:11}. A non-parametric approach based on Empirical Likelihood (EL) for testing a change in a linear model is considered by \cite{Liu:Zou:Zhang:08}. Always using the EL  method, the papers \cite{Zi:Zou:Liu:12}, \cite{Yu:Niu:Xu:13} construct the confidence region for the coefficient difference of a two-sample linear regression model. For a linear quantile model, \cite{Qu:08} proposes two types of statistics: one based on the subgradient and an another based on Wald statistic. For a generalized linear models, a method based on maximum of score statistics is used in \cite{Antoch} to test the change in the regression parameters.\\
\hh In this paper, we consider the change-point problem in a  general nonlinear model, by the EL method. Then, the  framework  of \cite{Liu:Zou:Zhang:08} is generalized. One of the major  difficulties  for  nonlinear model (beside the linear model approach) is that, for finding the test statistic,  the corresponding score functions  depend on the regression parameters, and above all, the analytical form of these derivatives is unknown. On the other hand, for linear models, many proofs are based on the convexity of the regression function with respect to the parameter regression, then, the extreme value of a convex function is attained on the boundary. These two factors lead to a more difficult theoretical study  of the test statistics for nonlinear model. Another difficulty to study the properties of the test statistic, for detecting a change in model, is due to the dependence on the change-points of the regression parameter estimator. To the authors' knowledge, the only paper which studies a hypothesis test in a change-point nonlinear  model is   \cite{Boldea:Hall:13} for very smooth nonlinear functions, using the least square method. But the least square method, in respect to the EL method, has the disadvantage that is  less efficient for outliers data. This occurs in the case of fatter tailed distributions of the error term.  Moreover, we will see in Section \ref{sec1:1} that the considered assumptions in \cite{Boldea:Hall:13} are stronger than in the present paper. \\
\hh Note that the paper  \cite{Hall:Sen:99} tests the structural stability in a nonlinear model by a generalized method of moments, and not a change in the regression parameters.\\
\hh I would emphasize that in the present paper, we have obtained an interesting result concerning the numerical simulations. The EL test outperforms the change detection by least square(LS) test proposed by \cite{Boldea:Hall:13}. The LS test does not work when the change-point is off-centred in the measurement interval. The proposed EL test does not this defect.  \\

The  paper is organized as follows. We first construct in Section \ref{sec1:1} a statistic, in order to test the change in the regression parameters of the  nonlinear model. The asymptotic behaviour of the test statistic  under the null hypothesis as well as under the  alternative hypothesis is studied. A particular case of two change-point model, the epidemic model, is considered in Section \ref{sec2:1}. In Section \ref{sec3:1}, simulations results illustrate the performance of the proposed test, concerning the empirical size, the empirical power and the estimation of the time of change, in particular when the error distribution is not Gaussian, when it has outliers or a large standard deviation. Some lemmas and their proofs are given in the last section (Appendix, Section \ref{sec4:1}).

\section{Test of a change-point }
\label{sec1:1}
\hh In this section, for a nonlinear model we are going to test the hypothesis that there is no change in the parameters of model (\ref{eqq1}) against the hypothesis that the parameters change from $\eb_1$ to $\eb_2$ at an unknown observation k, i.e. the model (\ref{eqq2}).
\subsection{Hypothesis, notations, assumptions}
\hh All throughout the paper, C denotes a positive generic constant which may take different values in different formula or even in different parts of the same formula. 
 All vectors are column and  $\textbf{v}^t$ denotes the transposed of $\textbf{v}$. All vectors and matrices are in bold.  Concerning the used norms, for a m-vector $\ev=(v_1, \cdots,v_m)$, let us denote by $\|\ev\|_1= \sum   _{j=1} ^ m |v_j|$ its $L_1$-norm and $\|\ev\|_2=(\sum   _{j=1} ^ m v^2_j)^{1/2}$ its $L_2$-norm. For a matrix $ \textbf{M}=(a_{ij})_{\substack{1\leqslant i \leqslant m_1\\1 \leqslant j \leqslant m_2}}$, we denote by $ \| \textbf{M}\|_1 = \max   _ {j=1,\cdots, m_2} (\sum  _{i=1}^{m_1} |a_{ij}|)$, the subordinate norm to the vector norm $\| .\|_1$.   Let $ \overset{\cal L} {\underset{n \rightarrow \infty}{\longrightarrow}}$, $ \overset{\eP} {\underset{n \rightarrow \infty}{\longrightarrow}}$, $ \overset{a.s.} {\underset{n \rightarrow \infty}{\longrightarrow}}$ represent convergence in distribution, in probability and almost sure, respectively, as $n \rightarrow \infty$.   \\
For coherence, we try to use the some notations as in the paper \cite{Liu:Zou:Zhang:08}, where the linear model was considered. This will allow to highlight the difficulties and results due to the nonlinearity. \\  \\
\hh  For each observation $i$, $Y_i$ denotes the response variable, $\XX_i$ is a $p \times 1$ random vector of regressors with distribution function $H(\xx)$, with $\xx \in \Upsilon  $, $ \Upsilon \subseteq \eR^p$, and $\varepsilon_i$ is the error. \ \\
\noindent The  continuous random vector sequence $(\XX_i,\varepsilon_i)_{1 \leq i \leq n}$ is independent identically distributed (i.i.d), with the same joint distribution as $(\XX,\varepsilon)$. For all $i$, $\varepsilon_i$ is independent of $\XX_i$.\ \\
The regression function $f:\Upsilon \times \Gamma \rightarrow \R$, with  $\Upsilon \subseteq \R^p$ and $\Gamma \subseteq \R^d$, is known up to a parameter $\eb=(\beta_1, \cdots,\beta_d)$. The parameter set $\Gamma$ is supposed compact. \ \\
In following, for $\xx \in \Upsilon$ and  $\eb \in \Gamma$, we use notation $\ef(\xx,\eb) \equiv \partial f(\xx,\eb)/\partial \eb$, $\eff(\xx, \eb)\equiv \partial^2 f(\xx,\eb)/\partial \eb^2$. \ \\  \\ 
\hh  With regard to the random variable $\varepsilon$ we make following assumption :\\
{\bf (A1)} $\eE[\varepsilon_i]=0$ and $\eE[\varepsilon_i^2] < \infty$, for all $i=1, \cdots,n$.\\
The regression function $f:\Upsilon \times \Gamma \rightarrow \R$ and the random vector $\XX$ satisfy the conditions :\\
\textbf{(A2)} for all $\xx \in \Upsilon$ and for $\eb \in \Gamma$, the function $f(\xx,\eb)$ is thrice differentiable in $\eb$ and continuous on $\Upsilon$.\\
\noindent \textbf{(A3)} $(| \frac{\partial ^2 f(\xx,\eb)}{\partial \beta_j \partial \beta_k }|)_{1 \leq j,k \leq d}$ and $(| \frac{\partial ^3 f(\xx,\eb)}{\partial \beta_j \partial \beta_k \partial \beta_l}|)_{1 \leq j,k,l \leq d}$ are bounded for any $\xx\in \Upsilon$ and $\eb$ in a neighbourhood of $\eb^0$. \\
\noindent  \textbf{(A4)} $\eE[\| \ef(\XX,\eb) \|_1]< \infty$, 
$\eE[\| \ef(\XX,\eb)\ef^t(\XX,\eb)\|_1]<\infty$ and $\eE[ | \frac{ \partial ^2 f(\XX,\eb) }{\partial \beta_j \partial \beta_k}| ]<\infty$, for all $ 1 \leq j,k \leq d$ and $\eb$ in a neighbourhood of $\eb^0$. \\

Assumptions (A3), (A4) are standard conditions, which are used in nonlinear models, for example see book \cite{seber:wild}. We remark that  assumption (A4) is weaker than the corresponding assumption employed in  paper \cite{Boldea:Hall:13}, where the least square method is used to test $H_0$ against $H_1$. The assumptions of the paper \cite{Boldea:Hall:13} are:  $\eE [ \varepsilon_i ^{2s}]$ $  <\infty$, $\,\, \sup  _ {\eb} \eE [ f (\XX,\eb)]^{2s} $  $< \infty$,  $ \,\, \sup   _ {\eb} \eE [ \| \ef (\XX,\eb)] \|_2] ^{2s} < \infty$ and for all $ 1 \leq j,k \leq d$,  $\,\, \sup _ {\eb} \eE [ | \frac{\partial ^2 f(\XX,\eb) }{\partial \beta_j \partial \beta_k}|^s] $ $< \infty$, for some $s>2$.\\
%
%$\sup  \limits _ {\eb} \eE [f (\xx,\eb)]^{2s} < \infty$,  $\sup  \limits _ {\eb} \eE [\| \ef (\xx,\eb)\|_2]^{2s} < \infty$,    $\sup   \limits_ {\xx, \eb} \eE [\eff (\xx,\eb)]^{2s} < \infty$, for some $s>2$.\\

We are interested in testing of the null hypothesis of no change in the model (\ref{eqq2}). Then the model has the form (\ref{eqq1}), that is 
\begin{equation}
\label{H0}
H_0 : \eb_1 =\eb_2=\eb.
\end{equation}
The alternative hypothesis assumes that one change occurs in the regression parameters, that is 
\begin{equation}
\label{H1}
H_1 : \eb_1 \neq \eb_2.
\end{equation}
Let $\ebo$ denote the true (unknown) of the parameter $\eb$ under hypothesis $H_0$ and $\eb^0_1$, $\eb^0_2$ (also unknown) the true parameters under hypothesis $H_1$.\\

In addition to the notations introduced above, let us consider the following $d$-random vectors 
 \begin{equation*}
\eg(\XX_i,\eb) \equiv \eg_i(\eb) \equiv \ef(\XX_i, \eb) [Y_i - f(\XX_i, \eb)].
\end{equation*}
We remark that, under the hypothesis $H_0$, we have $\eg_i(\ebo)=\ef(\XX_i,\ebo) \varepsilon_i$, for all $i=1, \cdots, n$ and $\eE[\eg_i(\ebo)]=0$. Consider also the $d \times d$ matrix
 \begin{equation*}
 \eV \equiv \eE[ \ef(\XX_i,\ebo) \ef^t(\XX_i,\ebo)].
\end{equation*}
Then $\Var(\varepsilon_i \ef(\XX_i,\ebo))=\sigma^2 \eV$.\\

 In order to introduce the empirical likelihood, let $y_1, \cdots, y_k, y_{k+1}, \cdots, y_n$ be observations for the random variables $Y_1, \cdots, Y_k,$ $Y_{k+1}, \cdots, Y_n$. For more details concerning empirical likelihood method, the reader can refer to \cite{owen}. Consider the following sets $I \equiv \{1, ..., k \}$ and $J\equiv \{ k+1, ..., n \}$, which contain the observation subscripts of the two segments for the model (\ref{eqq2}). Corresponding to these sets, let be the probability vectors $(p_1,\cdots,p_k)$ and $(q_{k+1},\cdots,q_{n})$. These vectors contained the probability to observe the value $y_i$ (respectively $y_j$) for the dependent variable $Y_i$ (respectively $Y_j$) : $p_i\equiv \eP [Y_i=y_i]$, for $i=1, \cdots, k$ and $q_j \equiv \eP[Y_j=y_j]$, for $j=k+1, \cdots,n$. Obviously, these probabilities satisfy the relations $\sum  _{i \in I}p_i = 1$ and $\sum  _{j \in J}q_j = 1$. \\

\subsection{Test statistics}
\hh  Under  hypothesis $H_0$ given by (\ref{H0}), the profile empirical likelihood (EL) for $\eb$ is   
 \begin{eqnarray*}
\begin{array}{ccl}
{\cal R}_{0,nk}(\eb)= \sup   _{(p_1,\cdots,p_k)}
  \sup   _{(q_{k+1},\cdots,q_n)}  \bigg\{  \prod    _{i \in I}p_i   \prod    _{j \in J}q_j ;   \sum    _{i\in I} p_i = 1,  \sum    _{j\in J} q_j =1, \\
  \sum    _{i \in I}p_i \eg_i( \eb )= \sum    _{j \in J}
q_j \eg_j( \eb)=\textbf{0}_d \bigg \},
\end{array}
\end{eqnarray*}
with $\textbf{0}_d$ the $d$-vector with all components zero. 
Without constraints $ \sum  _{i\in I} p_i \eg_i(\eb)=\textbf{0}_d$, the maximum of $\prod   _{i\in I} p_i$, $\prod   _{j \in J}q_j$ are  attained for $p_i=k^{-1}$, $q_j=(n-k)^{-1}$, respectively. Then, the profile EL ratio for $\eb$ has the form 
\begin{eqnarray}
\begin{array}{ccl}
\label{R'0}
{\cal R}'_{0,nk}(\eb)= \sup  _{(p_1,\cdots,p_k)}
  \sup    _{(q_{k+1},\cdots,q_n)} \bigg \{  \prod    _{i \in
I}k p_i   \prod  _{j \in J}(n-k) q_j ;   \sum    _{i\in I} p_i = 1, \\  \sum    _{j\in J} q_j =1, \sum    _{i \in I}p_i \eg_i( \eb )= \sum     _{j \in J}q_j \eg_j( \eb)=\textbf{0}_d  \bigg\}.
\end{array}
\end{eqnarray}
Similarly, under hypothesis $H_1$ given by (\ref{H1}) , the profile EL is 
\begin{eqnarray*}
\begin{array}{ccl}
{\cal R}_{1,nk}(\eb_1,\eb_2)= \sup    _{(p_1,\cdots,p_k)}
  \sup    _{(q_{k+1},\cdots,q_n)}  \bigg\{  \prod    _{i \in
I}p_i   \prod   _{j \in J}q_j ;   \sum    _{i\in I} p_i =1,  \\   \sum_{j\in J} q_j =1, \sum_{i \in I}p_i \eg_i( \eb_1 )= \textbf{0}_d, \sum_{j \in J}q_j \eg_j( \eb_2)=\textbf{0}_d  \bigg\}.
\end{array}
\end{eqnarray*}
Then, the profile EL ratio for $\eb_1, \eb_2$ has the form
\begin{eqnarray*}
\begin{array}{ccl} 
{\cal R}'_{1,nk}(\eb_1,\eb_2)= \sup   _{(p_1,\cdots,p_k)}\sup  _{(q_{k+1},\cdots,q_n)}
\bigg \{  \prod    _{i \in I}k p_i \prod   _{j \in J}(n-k) q_j
 ;   \sum   _{i\in I} p_i =1, \\  \sum   _{j\in J} q_j =1,  \sum_{i \in I}p_i \eg_i( \eb_1)=\textbf{0}_d,\sum   _{j\in J}q_j \eg_j( \eb_2)=\textbf{0}_d \bigg \}.
 \end{array}
\end{eqnarray*}
Thus, using an idea similar to the maximum likelihood test for testing $H_0$ against $H_1$, we consider the profile EL ratio 
\begin{equation}
\label{eq2}
\frac{{\cal R}_{0,nk}(\eb)}{{\cal R}_{1,nk}(\eb_1,\eb_2)}=\frac{{\cal R}'_{0,nk}(\eb)}{{\cal R}'_{1,nk}(\eb_1,\eb_2)},
\end{equation}
but, under this form, it has a complicated expression. 
In order to find a simpler form for the test statistic, we will study the denominator behaviour of the process given by (\ref{eq2}). \ \\

%********************************************************************
%********************************************************************%********************************************************************%********************************************************************%********************************************************************%********************************************************************%********************************************************************

The following result is a generalization of the nonparametric version of the Wilks theorem. More specifically, under $H_1$ due to the observation independence, on each segment we have a Wilks theorem. Then, we prove that, under $H_1$, the profile EL ratio for $\eb_1$, $\eb_2$ has a $\chi^2$ asymptotic distribution.
\begin{theorem}
\label{theo1} 
Suppose that assumptions (A1)-(A3) hold. Under the hypothesis $H_1$, we have 
$$ 
-2 \log {\cal R}'_{1,nk}(\eb_1,\eb_2) \overset{{\cal L}} {\underset{n \rightarrow \infty}{\longrightarrow}} \chi^2(2d).
$$
\end{theorem} 
\textit{Proof.} Under hypothesis $H_1$, on the first segment generated by the observations for $i \in I$, the profile EL function for $\eb_1$, for fixed $k$, is 
$$
\sup  _{(p_1,\cdots,p_k)} \bigg \{ \prod_{i \in
I}k p_i ; 0 \leq p_i \leq 1, \sum  _{i \in I} p_i=1,\sum _{i \in I}p_i \eg_i( \eb_1 )=\textbf{0}_d  \bigg \}.
$$
Applying the Lagrange multiplier method, using the paper \cite{Qin:law}, we have
$$
p_i=\frac{1}{k(1+\el^t_1 \eg_i(\eb_1))},
$$
where $\el_1 \in \eR ^ d$ the Lagrange multiplier. \\
Similarly, the profile EL function on the second segment generated by the observations for $j \in J$, is 
\begin{equation*}
\sup  _{(q_{k+1}, \cdots, q_n)} \bigg\{\prod  _{j \in
J}(n-k) q_j ; 0 \leq q_j \leq 1, \sum  _{j \in J} q_j=1,\sum  _{j \in J}q_j \eg_j( \eb_2 )=\textbf{0}_d \bigg \}.
\end{equation*}
\noindent This function is maximed for $q_j=(n-k-\el^t_2 \eg_i(\eb_2))^{-1}$, with $\el_2 \in \R^p$ the Lagrange multiplier. Then the empirical log-likelihood ratio statistic can be written  
\begin{equation}
\label{eq58}
 -2 \log {\cal R}'_{1,nk}(\eb_1,\eb_2)= 2\sum   _{i \in I}\log
\big[{1+ \el_1^t\eg_i(\eb_1 )}\big]+ 2\sum   _{j \in J}\log \big[{1-\el_2^t\eg_i( \eb_2 )}\big]  .
\end{equation} 
In view of Theorem 2 of \cite{Qin:law}, using assumptions (A1),
(A2) and (A3), each sum of the right-hand side of (\ref{eq58}) converges in law to $\chi^2(d)$. Taking into account that  the two terms of relation (\ref{eq58}) involved two independent sets of random vectors we obtain the theorem. \hspace*{\fill}$\blacksquare$  \ \\

%********************************************************************
%********************************************************************%********************************************************************%********************************************************************%********************************************************************%********************************************************************%********************************************************************

Consequently of this theorem, under hypothesis $H_1$ the denominator of the EL ratio given by (\ref{eq2}), is not asymptotically depend on the parameters $\eb_1$ and $\eb_2$. Then, from now on, we are going to consider that test statistic $-2 \log {\cal R}'_{0,nk}(\eb)$. \ \\

Taking into account the expression of ${\cal R}'_{0,nk}(\eb)$ given by  (\ref{R'0}) and using the Lagrange multiplier method, we have that  maximizing $ -2 \log {\cal R}'_{0,nk}(\eb)$ is equivalent to maximizing the following statistic with respect to $\eb$, $\eta_1$, $\eta_2$, $\el_1$, $\el_2$
\begin{equation}
\label{eqmu1}
\sum  _{i \in I}\big[ \log p_i - n \el_1^t p_i \eg_i(\eb)\big]+\sum  _{j \in J} \big[\log q_j +n \el^t_2  q_j \eg_j(\eb)\big]+\eta_1 \big(\sum  _{i \in I} p_i -1\big)+\eta_2 \big(\sum   _{j \in J} q_j -1\big),
\end{equation}
where $\eb \in \Gamma$, $\eta_1, \eta_2 \in \R$ and $\el_1, \el_2 \in \R^d$. \ \\
Since the derivatives of (\ref{eqmu1}) with respect to $p_i$, $q_j$ are null, using a similar argument as in the proof of Theorem \ref{theo1}, we obtain that
 \begin{equation}
\label{piqj}
p_i=\frac{1}{k+n \el^t_1 \eg_i(\eb)} \qquad q_j=\frac{1}{n-k- n \el^t_2 \eg_j(\eb)}.
\end{equation}
Then, the statistic $-2 \log {\cal R}'_{nk,0}(\eb)$ becomes 
\begin{equation}
\label{eq4} 
 2 \sum  _{i \in I}\log
\big[{1+\frac{n}{k}\el_1^t \eg_i( \eb)}\big]+ 2\sum  _{j \in J}\log \big[{1-\frac{n}{n-k}\el_2^t \eg_j(\eb )}\big].
\end{equation}
Taking into account relation (\ref{piqj}) for the probabilities $p_i$ and $q_j$, the derivative with respect to $\eb$ of (\ref{eq4}) is $2 n \big[ \sum  _{i \in I} p_i \el^t_i \egg_i(\eb)- \sum   _{j \in J}  q_j \el^t_2 \egg_j(\eb)\big]=0 $, with $\egg_i(\eb)$ the $d \times d$ matrix of the derivatives of vector $\eg_i(\eb)$ with respect to $\eb$, for $i=1,\cdots,k$. In order to have single parameters $\el$, we restrict the study to a particular case, when $\el_1$ and $\el_2$ satisfy the constraint $\eV_{1n}(\eb) \el_1 =\eV_{2n}(\eb) \el_2 $, with 
$$\eV_{1n}(\eb) \equiv k^{-1} \sum  _{i \in I} \egg_i(\eb),\qquad  \eV_{2n}(\eb) \equiv (n-k)^{-1} \sum   _{j \in J} \egg_j(\eb).$$
In the case of the true parameter $\ebo$, this two last matrices are denoted $\eV^0_{1n} \equiv \eV_{1n}(\ebo)$ and $\eV^0_{2n} \equiv \eV_{2n}(\ebo)$. Considering this constraint,  statistic (\ref{eq4}) becomes 
\begin{equation}
\label{eq5}
2 \sum_{i \in I}\log
\big[1+\frac{n}{k} \el^t\eg_i(\eb)\big]+ 2\sum  _{j \in J}\log
\big[1-\frac{n}{n-k}\el^t \eV_{1n}(\eb)\eV_{2n}^{-1}(\eb)\eg_j( \eb )\big].
\end{equation}
On the other hand, in order that the parameters belong a bounded set, in the place of $k$, we consider $\theta_{nk} \equiv k/n$, and we denote  statistic (\ref{eq5}) by $Z_{nk}(\theta_{nk},\el,\eb)$. Under hypothesis $H_1$, if $k_0$ is the point where the model change, we denote  $\theta_{0n}=k_0/n$. \ \\  \\
\hh Similar to the classical maximum likelihood test, but for models without change-points, we will study the maximum of empirical log-likelihood test statistic. For this, we calculate the score functions of test statistic (\ref{eq5})
\begin{eqnarray}
\label{eq6}
\ephi_{1n}(\theta_{nk}, \el, \eb) & \equiv & \frac{\partial Z_{nk}(\theta_{nk},\el,\eb)}{2 \partial \el}  \nonumber \\ &&= \sum  _{i \in I}
\frac{\eg_i(\eb)}{\theta_{nk}+\el^t \eg_i(\eb)}-\sum  _{j\in J}
\frac{\eV_{1n}(\eb)\eV_{2n}^{-1}(\eb)\eg_j(\eb)}{1-\theta_{nk}-\el^t\eV_{1n}(\eb)\eV_{2n}^{-1}(\eb)\eg_j(\eb)}.\nonumber \\
\end{eqnarray}
\begin{eqnarray}
\label{eq7}
 \ephi_{2n} (\theta_{nk},\el,\beta) & \equiv & \frac{\partial Z_{nk}(\theta_{nk},
\el,\eb)}{2\partial\eb}
 \nonumber \\ && = \sum  _{i\in I} \frac{
\egg_i(\eb) \el^t}{\theta_{nk}+
\el^t(\eb)\eg_i(\eb)}  
-\sum   _{j\in
J}  \frac{\partial (\eV_{1n}(\eb) \eV_{2n}^{-1}(\eb)\eg_j(\eb))/ \partial \eb }
{1-\theta_{nk}-\el^t \eV_{1n}(\eb) \eV_{2n}^{-1}(\eb)\eg_j(\eb)}\el^t. \nonumber \\
% \frac{\partial [\eV_1(\eb) \eV_2^{-1}(\eb)\eg_j(\eb)]}{\partial \eb} \nonumber
\end{eqnarray}
Then,  solving the system $\ephi_{1n}(\theta_{nk}, \el, \eb)=\textbf{0}_d $ and $\ephi_{2n}(\theta_{nk}, \el, \eb)=\textbf{0}_d$, the obtained solutions $\hat \el(\theta_{nk})$ and $\hat \eb(\theta_{nk})$ are the maximizers of the statistic (\ref{eq5}). We so obtain the profile maximum empirical likelihood function $Z_{nk}(\theta_{nk}, \hat \el(\theta_{nk}), \hat \eb(\theta_{nk}))$, which depends only on the change-point parameter $\theta_{nk}$. \ \\

 We emphasise that, compared with a linear model, in our case, matrix $\eV_{1n}(\eb)$, $\eV_{2n}(\eb)$ and derivative $\egg(\eb)$ depend on $\eb$. These, besides the nonlinearity of $\eg(\eb)$ involve difficulties in the study of the statistic $Z_{nk}(\theta_{nk},  \el,  \eb)$ and of the solutions $\hat \el(\theta_{nk})$, $\hat \eb(\theta_{nk})$.
\subsection{Asymptotic behaviour of the test statistic}

\hh In this section, for the probabilities given by (\ref{piqj}), under the constraint $\eV_{1n}(\eb) \el_1=\eV_{2n}(\eb) \el_2 $, we will first prove that $kp_i$, $(n-k)q_j$, can be framed by two strictly positive constants. This implies that the test statistic $Z_{nk}(\theta_{nk}, \hat \el(\theta_{nk}), \hat \eb(\theta_{nk}))$ is well defined. \ \\
 Properties established for $\hat \el(\theta_{nk})$ and $\hat \eb(\theta_{nk})$, solutions of (\ref{eq5}), will allow to consider instead of (\ref{eq5}), a more simple test statistic, given by relation (\ref{eqTn}). Next, we will study the asymptotic behaviour of this statistic, firstly under the hypothesis $H_0$ and next under $H_1$.
\subsubsection{Asymptotic behaviour under $H_0$}
\hh In order to study asymptotic behaviour of $\hat \el(\theta_{nk})$ and $\hat \eb(\theta_{nk})$, we will first study $kp_i$, for $i \in I$, and $(n-k)q_j$, for $j \in J$, with $p_i$, $q_j$ given by (\ref{piqj}). More exactly, we show that, if $\eb$ in the neighbourhood of $\ebo$, $kp_i$ and $(n-k)q_j$ can be framed by two strictly positive and bounded constants, with probability close to one.
\begin{proposition}
\label{prop}
Let the $\eta$-neighbourhood of $\ebo$, ${\cal V}_{\eta}(\ebo)= \{ \eb \in \Gamma; \| \eb-\ebo\|_2 \leq \eta \}$, with $\eta \rightarrow 0$. Under hypothesis $H_0$, suppose that assumptions (A1)-(A4) hold. Then we have\\
 (i) For all $i \in I$, for all $\epsilon>0$, there exist two constants $\textit{M}_1, \textit{M}_2>0$, such that, for all $\eb \in {\cal V} _{\eta}(\ebo)$,
\begin{equation}
\label{13}
 \eP \Big [ \frac{1}{M_2} \leq
\frac{1}{1+\frac{\el^t}{\theta_{nk}}g_i(\eb)}\leq  \frac{1}{M_1}
\Big ] \geq 1-\epsilon.
\end{equation}
(ii) For all $j \in J$, for all $\epsilon>0$,
there exist two constants  $\textit{M}_3, \textit{M}_4>0,$ such that, for all $\eb \in {\cal V} _{\eta}(\ebo)$,
\begin{equation}
\label{14}
 \eP \Big [ \frac{1}{M_4} \leq
\frac{1}{1-\frac{\el^t}{1-\theta_{nk}} \eV_{1n}(\eb) (\eV_{2n}(\eb))^{-1} 
\eg_j(\eb)}\leq \frac{1}{M_3} \Big ] \geq 1-\epsilon.
\end{equation}
\end{proposition}

%Now, we will prove the relations (\ref{13}) and (\ref{14}).

\noindent \textit{Proof. (i)} We consider the following decomposition for the Lagrange multiplier: $\el=\rho \ephi$, such that $\rho \geq 0$ and $ \|\ephi
\|_1=1$. Lemma \ref{lem3} implies that, there exists $M_2>0$, such that
$$\frac{1}{1+\frac{\el^t}{\theta_{nk}}\eg_i(\eb)}\geq \frac{1}{1+\frac{\rho}{\theta_{nk}}\|
\ephi^t \eg_i(\eb)\|_1} \geq \frac{1}{1+\frac{\rho}{\theta_{nk}}\|
\eg_i(\eb)\|_1} \geq
\frac{1}{M_2},$$
 with probability close to 1, that is, for all $\epsilon>0$,
\begin{equation}
\label{22}
 \eP[\frac{1}{1+ \frac{\rho}{\theta_{nk}}
\|\eg_i(\eb)\|_1} \geq \frac{1}{M_2}] \geq 1- \frac{\epsilon}{2}.
\end{equation}
For the right-hand side of relation (\ref{13}), we assume the contrary, that is, there exists $M_1>0$ such that 
$$\sup  _ {i \in I, \eb \in \Gamma}
 \frac{1}{1+\frac{\el^t}{\theta_{nk}}\eg_i(\eb)} \geq \frac{1}{M_1}.$$
This is equivalent to the fact that there exists $M_5>0$, such that 
\begin{equation*}
\inf  _ {i \in I, \eb \in \Gamma}
 \frac{\el^t}{\theta_{nk}}\eg_i(\eb) \leq - M_5.
\end{equation*}
 Since $\el=\rho \ephi$, $\rho>0$, and $0<\theta_{nk}<1$, therefore exists $M_6>0$ such that 
\begin{equation}
\label{eq23}
 \inf  _ {i \in I, \eb \in \Gamma}
  \ephi^t \eg_i(\eb)  \leq -M_6.
\end{equation} \ \\
On the other hand, we have that $ \inf   _ {i \in I, \eb \in \Gamma}\ephi^t \eg_i(\eb) \geq  -\inf   _ {i \in I, \eb \in \Gamma}\|  \eg_i(\eb) \|_1 $, with probability 1. Taking into account relation (\ref{eq23}), there exists $M_6>0$ such as $ -\inf  _ {i \in I, \eb \in \Gamma}\|  \eg_i(\eb) \|_1 \leq -M_6$ again too $  \sup   _ {i \in I, \eb \in \Gamma}\|  \eg_i(\eb) \|_1 \geq M_6$, which is in contradiction with relation ( \ref{eqg}). Then, the relation (\ref{13}) holds.\ \\
\textit{ (ii)} Relation (\ref{14}) can be proved in a similar way. \hspace*{\fill}$\blacksquare$  \ \\ 

%****************************************************************************
%****************************************************************************
%****************************************************************************
%****************************************************************************
%****************************************************************************
%****************************************************************************

By the following result, we show that $\hat{\el}(\theta_{nk})$ and $\hat{\eb}(\theta_{nk})$, the solutions of the score equations $\ephi_{1n}(\theta_{nk},\el,\eb)=\textbf{0}_d$ and $\ephi_{2n}(\theta_{nk},\el,\eb)=\textbf{0}_d$, have suitable properties. More precisely, we show that  $\|\hat{\el}(\theta_{nk})\|_2 \rightarrow 0$, as $n\rightarrow \infty$ and that $\hat{\eb}(\theta_{nk})$ is a consistent estimator of $\eb^0$, under hypothesis $H_0$. We also obtain their convergence rate. This will allow us to propose a simpler test statistics instead of $Z_{nk}(\theta_{nk},\el,\eb)$.

%****************************************************************************
%****************************************************************************
%****************************************************************************

\begin{theorem}
\label{theo2}
Suppose that the assumptions (A1)-(A4) hold. Under the hypothesis $H_0$, we have
$\hat{\el}(\theta_{nk})=$ min $\{ {\theta_{nk}, 1-\theta_{nk}}\} O_{\eP}((n$ min $\{ {\theta_{nk}, 1-\theta_{nk}}\})^{-1/2})$  and $\hat{\eb}(\theta_{nk})-\eb^0=O_{\eP}((n$ min $\{ {\theta_{nk}, 1-\theta_{nk}}\})^{-1/2})$. 
\end{theorem}
\noindent {\it Proof.}
The structure of the proof is similar to that of linear model (Lemma A1 of \cite{Liu:Zou:Zhang:08}) but important modifications and supplementary results are necessary, due to the model nonlinearity. Without loss of generality, we assume that 
 $\min  \{\theta_{nk}, 1-\theta_{nk}\}= \theta_{nk}$. The other case is similar.\\
  By the definition of the profile empirical likelihood ratio ${\cal R}'_0(\eb)$, we have the following constraints 
\begin{equation}
\label{eq9}
\textbf{0}_d=\sum  _{i \in I} p_i \eg_i(\eb)=\sum  _{j \in J} q_j \eg_j(\eb).
\end{equation}
 We recall that, under hypothesis $H_0$, the expression of $p_i$ is given by (\ref{piqj}), and it is equal to $(\theta_{nk}+n \el^t \eg_i(\eb))^{-1}$, for $i =1, \cdots, n\theta_{nk}$. 
Then, by elementary calculations, we obtain 
\begin{equation}
\label{eq10}
\textbf{0}_d =\frac{1}{n \theta_{nk}} \sum  _{i \in I} \eg_i(\eb) -\frac{1}{n \theta_{nk}^2} \sum  _{i \in I} \frac{\eg_i(\eb) \eg_i^t (\eb) }{1+\frac{\el^t(\eb)}{\theta_{nk}}\eg_i(\eb)}\el(\eb). 
\end{equation}
Let us make the remark that we denote $\el$ by $\el(\eb)$ in order to indicate that for each value of $\eb$, solution of (\ref{eq10}), we will have a different value for $\el$. 
We take $\eb=\eb^0 \pm(n \theta_{nk})^{-r}\textbf{1}_d$, with $\textbf{1} _d $ the $d$-vector with all components 1 and $r>0$ will be specified later. Therefore, $\| \eb-\ebo \|_2 =(n \theta_{nk})^{-r} \rightarrow 0 $, as $n \theta_{nk} \rightarrow \infty$.\ \\ 
\hspace*{0.4cm} For the first sum of the right-hand side of (\ref{eq10}), by Lemma \ref{lem4}, we have 
\begin{equation*}
\frac{1}{n \theta_{nk}} \sum  _{i \in I} \eg_i(\eb)=O_{\eP}((n \theta_{nk})^{-1/2})+\eV_{1n}^0 (\eb-\ebo)+o_{\eP}(\eb-\eb^0).
\end{equation*}

 Now, we consider the second term of the right-hand side of relation (\ref{eq10}). From Proposition \ref{prop}, we have that for all $\epsilon>0$, there exists $M_1, M_2 >0$, such that 
\begin{equation*}
\label{eq13}
\eP \Big [\frac{1}{M_1}  \sum _{i \in I} \eg_i(\eb) \eg^t_i(\eb)\leq  \sum   _{i \in I}\frac{\eg_i(\eb) \eg^t_i(\eb)}{1+\frac{\el^t_(\eb)}{\theta_{nk}}
 \eg_i(\eb)} \leq \frac{1}{M_2}  \sum  _{i \in I} \eg_i(\eb) \eg^t_i(\eb) \Big ]<\epsilon.
\end{equation*} 
 This implies that, in order to study the second term of the right-hand side of the relation (\ref{eq10}), we must study only $(n\theta_{nk})^{-1} \sum  _{i \in I} \eg_i(\eb) \eg^t_i(\eb)$. By a Taylor's expansion of $\eg_i(\eb)$ in a neighbourhood of $\ebo$, using an argument similar to the one used for the first term of (\ref{eq10}), together with the assumption (A3), we obtain 
\begin{equation}
\label{eq14}
\frac{1}{n\theta_{nk}} \sum  _{i \in I} \eg_i(\eb) \eg^t_i(\eb)=\frac{1}{n\theta_{nk}} \sum  _{i \in I} \eg_i(\ebo) \eg^t_i(\ebo)(1+o_{\eP}(1)).
\end{equation} 
Taking into account Lemma \ref{lem4} and relation  (\ref{eq14}), the relation (\ref{eq10}) becomes 
\begin{equation}
\label{eq15}
\textbf{0}_d= \Big  [ O_{\eP}((n \theta_{nk})^{-1/2})+ \eV_{1n}^0 (\eb-\ebo) -\frac{1}{n\theta_{nk}^2} \sum  ^{n \theta_{nk}}_{i=1} \eg_i(\ebo) \eg^t_i(\ebo)\el(\eb)\Big](1+o_{\eP}(1))
\end{equation}
 We consider a constant $r$ such that $1/3 \leq r <1/2$. If $\eb=\ebo +(n \theta_{nk}) ^{-r} \textbf{1}_d$, then
$(\eb-\ebo)^t \textbf{1}_d >0$, and if $\eb=\ebo - (n \theta_{nk})^{-r}\textbf{1}_d$ then $(\eb-\ebo)^t\textbf{1}_d<0$. Then, the relation (\ref{eq15}) implies 
\begin{eqnarray}
\label{eq16}
 \el(\ebo \pm (n \theta_{nk})^{-r} \textbf{1}_d)&=& \pm \Big [ \theta_{nk}   \Big( \frac{1}{n \theta_{nk}}
\sum   _{i \in I}\varepsilon_i^2 \ef_i( \ebo) \ef^t_i(\ebo) \Big)^{-1}  \eV_{1n}^0(n   \theta_{nk})^{-r} \textbf{1}_d \nonumber \\ && +O_{\eP}((n \theta_{nk})^{-1/2}) \Big ](1+o_{\eP}(1)).
\end{eqnarray}
 For the observations $j \in J$, let us consider the function $\ev : \Gamma \rightarrow \eR^d$ defined by 
\begin{equation*}
 \ev(\eb)=  \sum  _{j \in J} q_j \eg_j(\eb)=\frac {1}{n-n \theta_{nk}} \sum  _{j \in J}\frac
{\eg_j(\eb)}{1-\frac{\el^t(\eb)}{1-\theta_{nk}}\eV_{1n}(\eb) \eV_{2n}^{-1}(\eb) \eg_j(\eb)}.
\end{equation*}
 Note that $\ev(\hat \eb(\theta_{nk}))=\textbf{0}_d$. For $\ev(\eb)$, we have the following decomposition 
\begin{equation*}
\frac{ \eV_{1n}(\eb) \eV_{2n}^{-1}(\eb)}{n(1- \theta_{nk})^2 } \sum  _{j \in J} \frac{\eg_j(\eb) \eg^t_j(\eb)}{1-\frac{\el^t(\eb)}{1-\theta_{nk}}\eV_{1n}(\eb) \eV_{2n}^{-1}(\eb) \eg_j(\eb)}\el(\eb) +\frac{1}{n(1-\theta_{nk})} \sum  _{j\in J}  \eg_j(\eb). 
\end{equation*}
To facilitate writing, we consider the following $d \times d$ squares matrices, defined by 
\begin{equation}
\label{eqD1D2} 
\eD^0_{1n}=\frac{1}{n\theta_{nk}}\sum  _{i \in I} \eg_i(\ebo)\eg^t_i(\ebo), \,\,\,\,
\eD^0_{2n}=\frac{1}{n-n\theta_{nk}}\sum  _{j \in J} \eg_j(\ebo)\eg^t_j(\ebo). 
\end{equation}
As for the observations $i \in I$, we obtain, similarly as for relation (\ref{eq15}),   \\
 $\ev(\eb)=[\eV_{2n}^0 (\eb-\ebo)+ (1-\theta_{nk})^{-1} \eV_{1n}^0 (\eV_{2n}^0)^{-1}\eD^0_{2n} \el(\eb) +O_{\eP}((n(1-\theta_{nk}))^{-1/2})]$\\ 
 $\cdot (1+o_{\eP}(1)).$ \\ 
Replacing $\el(\eb)$ by the value obtained in (\ref{eq16}), we obtain \\
$ \ev(\eb)= [\eV_{2n}^0 ( \eb-\ebo) + (\theta_{nk})(1-\theta_{nk})^{-1}\eV_{1n}^0
  (\eV_{2n}^0)^{-1}
\eD_{2n}^0 (\eD_{1n}^0)^{-1} \eV_{1n}^0  ( \eb-\ebo)+O_{\eP}((n(1-\theta_{nk}))^{-1/2})+O_{\eP}((n \theta_{nk})^{-1/2}) ](1+o_{\eP}(1))$. \\
Because $\eb= \ebo \pm (n \theta_{nk})^{-r} \textbf{1}_d$, $1/3 \leq r <1/2$ and $\min\{ \theta_{nk}, 1- \theta_{nk}\}=\theta_{nk}$, then $v(\eb)$  becomes 
\begin{equation}
\label{eq18} 
 [( \eV_{2n}^0 + \frac{\theta_{nk}}{1-\theta_{nk}} \eV_{1n}^0
  (\eV_{2n}^0)^{-1}
\eD_{2n}^0 (\eD_{1n}^0)^{-1} \eV_{1n}^0  )( \eb-\ebo)+O_{\eP}((n \theta_{nk})^{-1/2})) ](1+o_{\eP}(1)).
\end{equation}
This implies that $\ev(\ebo + (n \theta_{nk})^{-r} \textbf{1}_d)$ and $\ev(\ebo-(n \theta_{nk})^{-r} \textbf{1}_d)$ have a different signs, component by component. Moreover, because $\ev$ contains continuous functions in the neighbourhood of $\ebo$, there exists a $\eb$ such that $\ev(\eb)=\textbf{0}_d$. But since $\ev(\hat{\eb}(\theta_{nk}))=\textbf{0}_d$, we have that $\hat{\eb}(\theta_{nk}) \in [\ebo- (n \theta_{nk})^{-r}\textbf{1}_d,\ebo+ (n \theta_{nk})^{-r} \textbf{1}_d]$, which implies, because $r<1/2$, that  $\hat{\eb}(\theta_{nk})-\ebo=O_{\eP}((n \theta_{nk})^{-r}) \geq O_{\eP} ((n \theta_{nk})^{-1/2})$. This last relation, together with the relation (\ref{eq18}), since $\hat \eb(\theta_{nk})-\ebo$ is the coefficient of a matrix strictly positive, implies that in order to have $\ev(\hat \eb(\theta_{nk}))=\textbf{0}_d$, we must have $\hat{\eb}(\theta_{nk})-\ebo=O_{\eP}((n \theta_{nk})^{-1/2})$. Considering this result, for the relation (\ref{eq16}), we obtain $\el(\hat \eb(\theta_{nk}))=\theta_{nk} O_{\eP}((n \theta_{nk})^{-1/2})$. The theorem is completely proved.\hspace*{\fill}$\blacksquare$ \\

%********************************************************************
%********************************************************************%********************************************************************%********************************************************************%********************************************************************%********************************************************************%********************************************************************

\begin{Remark}
\label{rem1}
In view of the proof of Theorem \ref{theo2}, under hypothesis $H_0$, we can consider instead of $Z_{nk}(\theta_{nk},\el,\eb)$, given by (\ref{eq5}), the following modified statistic 
\begin{equation}
\label{eqTn}
T_{nk}(\theta_{nk},\el,\eb)=2 \sum_{i \in I}\log
(1+\frac{1}{\theta_{nk}} \el^t\eg_i(\eb))+ 2\sum  _{j \in J}\log
(1-\frac{1}{1-\theta_{nk}}\el^t \eg_j( \eb )).
\end{equation}
\end{Remark} 

Because the regression function is nonlinear, in order to the maximum empirical likelihood always exists, we consider that the parameter $\theta_{nk} \in [\Theta_{1n}, \Theta_{2n}] \subset (0,1)$, such that  $n \Theta_{1n} \rightarrow \infty$, $n(1- \Theta_{2n}) \rightarrow \infty$, as $n \rightarrow \infty$ for example. The reader can find a discussion concerning the possible values of $\Theta_{1n}$, $\Theta_{2n}$ in the papers \cite{Zou:Liu:Qin:Wang:14}, \cite{Liu:Zou:Zhang:08}. \\
Finally, the test statistic  for testing the hypothesis $H_0$ against $H_1$ is 
\begin{equation}
\label{eq8}
\tilde T_n \equiv \max_{\theta_{nk} \in [\Theta_{1n}, \Theta_{2n}]} T_{nk}(\theta_{nk}, \hat \el(\theta_{nk}), \hat \eb(\theta_{nk})).
\end{equation}

Then, we can consider as estimator  for the time of change $k^0$, the maximum empirical likelihood estimator: $\tilde{k}_n\equiv n\tilde{\theta}_n\equiv n \min \{ \tilde{\theta}_{nk}; \tilde{\theta}_{nk}= \argmax   _{\theta_{nk} \in [\Theta_{1n}, \Theta_{2n}]}$ \\ $ T_{nk}(\theta_{nk}, \hat \el(\theta_{nk}), \hat \eb(\theta_{nk})) \}$. Recall that $\hat {\el}(\theta_{nk})$ and $\hat {\eb}(\theta_{nk})$ are the solutions of the score equations (\ref{eq6}) and (\ref{eq7}). \ \\

The following result gives the asymptotic distribution of the test statistic $\tilde {T}_n$ given by (\ref{eq8}), under the null hypothesis of no-change. For this purpose, we consider functions: $A(x) \equiv (2\log x)^{1/2}$,
$D(x)=2 \log x + \log \log x$ and $u(n)= \frac{1- \Theta_{1n} \Theta_{2n}}{\Theta_{1n}(1- \Theta_{2n})} \rightarrow \infty$ as $n  \rightarrow \infty$.   
%********************************************************************
%********************************************************************%********************************************************************%********************************************************************%********************************************************************%********************************************************************%********************************************************************

\begin{theorem}
\label{theo3} 
Under the assumptions (A1)-(A4), if the hypothesis $H_0$ is true, then we have, for all $t \in \R$ 
\begin{equation}
\label{eqth1}
\lim_{n\rightarrow \infty} \eP \{ A(\log
u(n))(\tilde{T}_n)^{\frac{1}{2}} \leq t+ D(\log u(n))\}= exp(-e^{-t}).
\end{equation}
\end{theorem}
\noindent {\it Proof.} The structure of the proof is similar to that of Theorem 1.3.1 of \cite{Csorgo:Horvath:97}, with the modification that the Lemma \ref{lem2} is used instead of Theorem 1.1.1 of \cite{Csorgo:Horvath:97}, and the Theorem A.3.4 of \cite{Csorgo:Horvath:97} instead of Corollary A.3.1 of \cite{Csorgo:Horvath:97}. The details are omitted. \hspace*{\fill}$\blacksquare$ \\ 
\begin{Corollary}
\label{corol1}
Consequence of this theorem, for a fixed size $\alpha \in (0,1)$, we can deduct the critical test region :
\begin{equation*}
(\tilde{T}_n)^{1/2} \geq \frac{-\log (-\log \alpha)+D(\log u(n))}{A(\log u(n))}.
\end{equation*}
\end{Corollary}
%**************************************************************************
%%%%%%%%%%%%%%%%%%%%%%%%%%%%%%%%%%%%%%%%%%%%%%%%%%%%%%%%%%%%%%%%%%%%%%%%%%%%%%%%
%**************************************************************************
%%%%%%%%%%%%%%%%%%%%%%%%%%%%%%%%%%%%%%%%%%%%%%%%%%%%%%%%%%%%%%%%%%%%%%%%%%%%%%%%
\vskip 0.5cm
\hh Using Theorem \ref{theo3} in applications is quite complicated. First, because we must first solve equation system (\ref{eq6}) and (\ref{eq7}) where the nonlinearity in parameter $\eb$ up to and including in matrices  $\eV_{1n}(\eb)$, $\eV_{2n}(\eb)$, $\eV_{2n}^{-1}(\eb)$ causes numerical difficulties and long computation time. Moreover, then it must then find $\theta_{nk}$ that maximizes statistic (\ref{eq8}). We can propose an approached form for the test statistic much simpler to use in practice, but which preserves the theoretical properties of (\ref{eq8}). 
\begin{Remark}
\label{rem2}
Taking into account the last relation of Lemma \ref{lem2}, Theorem \ref{theo3} implies that, in practice, for testing the hypothesis $H_0$ against $H_1$, we will use an approximate form 
\begin{equation}
\label{eqdez}
 T(\theta_{nk})= n
\sigma^{-2} \theta_{nk}(1-\theta_{nk})
(\textbf{W}_{1n}^0-\textbf{W}_{2n}^0)^t \eV^{-1}
(\textbf{W}_{1n}^0-\textbf{W}_{2n}^0)  
(1+o_{\eP}(1) ),
\end{equation} 
\noindent where 
\begin{equation}
\label{w1w2}
\textbf{W}_{1n}^0 =
\frac{1}{n\theta_{nk}} \sum  _{i \in I}\eg_i( \ebo), \qquad  \textbf{W}_{2n} ^0=
\frac{1}{n(1-\theta_{nk})}\sum  _{j \in J}\eg_j( \ebo). 
\end{equation}
Since $\ebo$ is unknown, in applications, we replace it with a consistent estimator, for example, the ordinary least square estimator, denoted by $\hat \eb _{LS}$. Under $H_0$, error variance $\sigma^2$ is estimated by $n^{-1} \sum   _{i=1}^n [\textbf{Y}_i - f(\XX_i,\hat \eb _{LS})]^2$ and matrix $\eV$ by $n^{-1} \sum   _{i=1} ^ n  \ef(\XX_i,\hat \eb _{LS})\ef^t(\XX_i,\hat \eb _{LS})$. \\
The approached maximum empirical likelihood estimator for the time of change $k^0$ is $$\hat{k}_n=n\hat{\theta}_n=n \min \{ \hat{\theta}_{nk}; \hat{\theta}_{nk}= \argmax   _{\Theta_{1n} \leq \theta_{nk} \leq \Theta_{2n}} T(\theta_{nk})\}.$$
\end{Remark}

%********************************************************************
%********************************************************************%********************************************************************%********************************************************************%********************************************************************%********************************************************************%********************************************************************
\subsubsection{Asymptotic behaviour of $T_{nk}$ and $\tilde T_{nk}$ under $H_1$}

\hh We consider now that the hypothesis $H_1$ is true. We will first prove that the maximum of statistic $T_{nk}$ converges to the maximum of its limit distribution. Then, we will show that statistic test $\tilde T_{n}$ is consistent (it has asymptotic power equal to 1). If $k^0$ is the true time of change, we denote by $\theta_{n0} =k^0/n$ and we suppose that $\theta_0 \equiv \lim \limits _{n \rightarrow \infty} \theta_{n0}$, where $\theta_0 \in (0,1)$.\\

 For $\xx \in \Upsilon$ and $e \in \eR$, let $F_{\xx}(e)$ and $G_{\xx}(e)$ the conditional distributions of $\eg(\XX_i,\eb)$ when $\XX_i=\xx$ for $i \in I$ and $j \in J$, respectively. Let $\e1_{(.)}$ the indicator function. Recall that, the distribution function of $\XX$ is $H(\xx)$. For $\xx$ and $\theta$ fixed, we define  \ \\ 
$dP_{\xx}{(e)} \equiv(\theta \e1_{\{\theta \leq \theta_0\}}+ \theta_0
\e1_{\{\theta
>\theta_0\}})dF_{\xx}(e) +(\theta - \theta_0) \e1_{\{\theta > \theta_0\}}
dG_{\xx}(e)$, \ \\
$dQ_{\xx}{(e)}\equiv ((1-\theta) \e1_{\{\theta \geq \theta_0\}}+
(1-\theta_0) \e1_{\{\theta < \theta_0\}})dG_{\xx}(e) + (\theta_0 -
\theta)\e1_{\{\theta < \theta_0\}} dF_{\xx}(e)$, \ \\
$dR_{\xx}{(e)}\equiv \e1_{\{\theta < \theta_0\}} dF_{\xx}(e) + \e1_{\{\theta > \theta_0\}} dG_{\xx}(e)$. \ \\ 

Since under $H_0$, we proved that instead of EL statistic (\ref{eq5}) we can consider statistic (\ref{eqTn}), let us define  the following statistic
\begin{equation}
\label{eqlam}
\Lambda_{nk} (\theta_{nk})=  T_{nk}(\theta_{nk}, \tilde{\el}(\theta_{nk}),
\tilde{\eb}(\theta_{nk}))/ (2n), \Lambda_n(0)= \Lambda_n(1)=0,
\end{equation}
with $T_{nk}$ given by relation (\ref{eqTn}), and $ \tilde{\el}(\theta_{nk}),
\tilde{\eb}(\theta_{nk})$ solutions of the system  
\begin{equation}
\label{derTn}
\left\{
\begin{array}{ccl}
\frac{\partial T_{nk}(\theta_{nk},\el,\eb)}{2 \partial \el}=
\sum  _{i \in I}
\frac{\eg_i(\eb)}{\theta_{nk}+\el^t \eg_i(\eb)}-\sum  _{j\in J}
\frac{\eg_j(\eb)}{1-\theta_{nk}-\el^t\eg_j(\eb)}=\textbf{0}_d, \\\frac{\partial T_{nk}(\theta_{nk},\el,\eb)}{2 \partial \eb}=
\sum  _{i \in I}
\frac{\egg_i(\eb) \el}{\theta_{nk}+\el^t \eg_i(\eb)}-\sum  _{j\in J}
\frac{\egg_j(\eb) \el}{1-\theta_{nk}-\el^t \eg_j(\eb)}=\textbf{0}_d.
\end{array}
\right.
\end{equation}
By a similar proof to that of Theorem \ref{theo2}, under $H_0$, we have that 
\begin{eqnarray}
\label{lbtild}
&& \tilde{\el}(\theta_{nk})= min \{ {\theta_{nk}, 1-\theta_{nk}}\} O_{\eP}((n \, min \{ {\theta_{nk}, 1-\theta_{nk}}\})^{-1/2}),
 \nonumber \\ && 
 \tilde{\eb}(\theta_{nk})-\eb^0=O_{\eP}((n \,min \{ {\theta_{nk}, 1-\theta_{nk}}\})^{-1/2}).
\end{eqnarray}
\hh For any $\el$ and $\eb$, let the function $K :  \Upsilon \times \eR \times (0,1)$ defined by 
$$K(\xx,e,\theta)= \theta + \el^t \ef(\xx,\eb)[e-f(\xx,\eb)+f(\xx,\ebo)].$$
Let  also 
\begin{eqnarray}
\psi (\theta,\el,\eb) & =& \int  _ {\Upsilon} \Big( \int  _ {\eR} \log
K(\xx,e,\theta) dP_{\xx}{(e)} + \int  _ {\eR} \log (1-K(\xx,e,\theta))
 dQ_{\xx}{(e)}\Big) dH(\xx) \nonumber \\ && - \theta \log \theta - (1-\theta) \log(1-\theta).
\end{eqnarray}
We will prove by Theorem \ref{theo4} that $\psi$ is the limit process of $\Lambda_{nk}$, under $H_1$.\\
Then consider, for a fixed $\theta \in (0,1)$, $ \tilde{\tilde{\el}}(\theta),
\tilde{\tilde{\eb}}(\theta)$  the solutions to the following score equations 
 \begin{equation}
\label{eqz1}
\left\{
\begin{array}{ccl}
\ez_1(\theta,\el,\eb)=  \int  _ {\Upsilon} \Big( \int  _ {\eR} \frac
{\eg(\xx,\eb)}{K(\xx,e,\theta)} dP_{\xx}{(e)} -  \int  _ {\eR} \frac {\eg(\xx,\eb)}{1-K(\xx,e,\theta)} dQ_{\xx}{(e)}
\Big) dH(\xx)= \textbf{0}_d, \\
\ez_2(\theta,\el,\eb)= \int  _ {\Upsilon} \Big( \int  _ {\eR} \frac {\egg(\xx,\eb)\el}{K(\xx,e,\theta)}
dP_{\xx}{(e)} - \int  _ {\eR} \frac {\egg(\xx,\eb)\el}{1-K(\xx,e,\theta)} dQ_{\xx}{(e)} \Big) dH(\xx)=\textbf{0}_d,
\end{array}
\right.
\end{equation}
where, $\ez_1(\theta,\el,\eb)=\partial \psi
(\theta,\el,\eb) /\partial \el$ and $\ez_2(\theta,\el,\eb)= \partial \psi (\theta,\el,\eb)/\partial \eb$. \ \\

 We require the following assumptions for the next theorems : \ \\

\noindent \textbf{(A5)} The matrix $ 
\begin{pmatrix}
   -\frac{\partial \ez_1(\theta,\el,\eb)}{\partial \el} & -\frac{\partial \ez_1(\theta,\el,\eb)}{\partial \eb} \\
 -\frac{\partial \ez_2(\theta,\el,\eb)}{\partial \el} & -\frac{\partial \ez_2(\theta,\el,\eb)}{\partial
\beta}
\end{pmatrix}
$ is nonsingular for all $\theta \in (0,1)$. 

%\noindent  \textbf{(A6)} \ \\
%\medskip $ \qquad \qquad \int  _ {\Upsilon} \int  _ {\eR} (\frac{\eg(\xx,\eb) \eg^t(\xx,\eb)}{K^2}+\frac{\eg(\xx,\eb) \eg^t(\xx,\eb)}{(1-K)^2})d(F(\xx,e)+G(\xx,e)) < \infty,$ \ \\
%\medskip $\qquad \qquad \int  _ {\Upsilon} \int  _ {\eR} ( \frac{\egg(\xx,\eb) \egg(\xx,\eb)}{K^2}+\frac{\egg(\xx,\eb) \egg(\xx,\eb) }{(1-K)^2})d(F(\xx,e)+G(\xx,e)) < \infty, \,\,\,\, \hbox{and}$  \ \\
%\medskip $ \qquad \qquad  \int  _ {\Upsilon}  \int  _ {\eR} (\frac{\eggg(\xx,\eb) }{K^2}+\frac{\eggg(\xx,\eb) }{(1-K)^2})d(F(\xx,e)+G(\xx,e)) < \infty$. \ \\

\noindent  \textbf{(A6)} 
The two following integrals are applied component by component of the corresponding matrix. \ \\
%Each term of the matrices $\eg(\xx,\eb) \eg^t(\xx,\eb)$ and $\egg(\xx,\eb) \egg(\xx,\eb)$ verify respectively, \ \\
\medskip $ \qquad \qquad \int  _ {\Upsilon} \int  _ {\eR} (\frac{\eg(\xx,\eb) \eg^t(\xx,\eb)}{K^2(\xx,e,\theta)}+\frac{\eg(\xx,\eb) \eg^t(\xx,\eb)}{(1-K(\xx,e,\theta))^2})d(F_{\xx}(e)+G_{\xx}(e)) < \infty,$ \ \\
\medskip $\qquad \qquad \int  _ {\Upsilon} \int  _ {\eR} ( \frac{\egg(\xx,\eb) \egg(\xx,\eb)}{K^2(\xx,e,\theta)}+\frac{\egg(\xx,\eb) \egg(\xx,\eb) }{(1-K(\xx,e,\theta))^2})d(F_{\xx}(e)+G_{\xx}(e)) < \infty.$ \ \\ % \,\,\,\, \hbox{and}  \ \\
Let $\eg_l(\xx,\eb)$ the l-th component of the vector $\eg(\xx,\eb)$. For all $ 1 \leq j, k \leq d$, we suppose that \ \\
\medskip $ \qquad  \int  _ {\Upsilon}  \int  _ {\eR} (\frac{1 }{K^2(\xx,e,\theta)}  \frac{\partial ^2 \eg_l(\xx,\eb) }{\partial \eb_j\partial \eb_k }+ \frac{1 }{(1-K(\xx,e,\theta))^2}  \frac{\partial ^2  \eg_l(\xx,\eb) }{\partial \eb_j\partial \eb_k })d(F_{\xx}(e)+G_{\xx}(e)) < \infty$. \ \\ \\
\noindent \textbf{(A7)} The functions $f(\xx,\eb)$ and $\ef(\xx,\eb)$ are equicontinuous in $\eb$ on $\Gamma$.  \\

\begin{Remark}
\label{rem3}
A sufficient condition for the equicontinuity of the functions $f(\xx,\eb)$ and $\ef(\xx,\eb)$ is that they are Lipschitzian with respect to $\eb$ on $\Gamma$.
\end{Remark}

%********************************************************************
%********************************************************************%********************************************************************%********************************************************************%********************************************************************%********************************************************************%********************************************************************

Following theorem shows that if $\theta_{nk}$ converges to the true value $\theta_0$, then the maximum of the modified EL test statistic converges to the maximum of its  limit distribution.

\begin{theorem}
\label{theo4}
 \noindent Under the alternative hypothesis $H_1$, if the assumptions (A5)-(A7) are satisfied, $\theta_0 \in (0,1)$ and $
 \lim  _{n\rightarrow \infty} \theta_{nk}= \theta
\in [0,1]$, then $\Lambda_{nk}(\theta_{nk})
\overset{{a.s.}} {\underset{n \rightarrow \infty}{\longrightarrow}} \psi( \theta, \tilde{\tilde{\el}}(\theta),
\tilde{\tilde{\eb}}(\theta))$, where $\psi( \theta, \tilde{\tilde{\el}}(\theta),
\tilde{\tilde{\eb}}(\theta))$ is
a strictly increasing function on $(0, \theta_0)$
 decreasing on $( \theta_0,1)$ and $\max  _
{0 \leq \theta \leq 1} \psi( \theta, \tilde{\tilde{\el}}(\theta)$,  $\tilde{\tilde{\eb}}(\theta))=\psi( \theta_0, \tilde{\tilde{\el}}(\theta_0),
\tilde{\tilde{\eb}}(\theta_0))$ \\ $=0$.
\end{theorem}

\noindent {\it Proof.} We will prove this theorem in three steps.\\
{\it Step  1.} We first prove that, for all fixed $\theta \in (0,1)$, we have 
\begin{equation}
\label{argmax1}
\argmax   _{(\el,\eb)} T_{nk} (\theta,\el,\eb)\overset{{a.s.}} {\underset{n \rightarrow \infty}{\longrightarrow}} \argmax   _{(\el,\eb)} \psi (\theta,\el,\eb). 
\end{equation}
Obviously, by the law of large numbers, for all $(\theta,\el,\eb) \in (0,1) \times \eR \times \Gamma$, we have $(2n)^{-1}  T_{nk} (\theta,\el,\eb)\overset{{a.s.}} {\underset{n \rightarrow \infty}{\longrightarrow}} \psi (\theta,\el,\eb) $. On the other hand, by the assumption (A5), $\argmax   _{(\el,\eb)} \psi (\theta,\el,\eb)$ is the unique solution of the system (\ref{eqz1}). Seen the assumptions (A6) and (A7), the function  
$(2n)^{-1}  T_{nk} (\theta,\el,\eb)$ is equicontinuous and bounded in $\el$ and $\eb$. Then, using Theorem 1.12.1 of \cite{Van Vaart:Wellner:25}, we have that the convergence of $(2n)^{-1}  T_{nk} (\theta,\el,\eb)$ to $\psi (\theta,\el,\eb)$ is uniform in $(\el,\eb)$. Taking into account that the solution of system (\ref{derTn}) is unique, we obtain relation (\ref{argmax1}). \ \\
{\it Step  2.} We show that 
\begin{equation}
\label{argmax2}
\max  _ {\theta_{nk}} \Lambda_{nk} (\theta_{nk})\overset{{a.s.}} {\underset{n \rightarrow \infty}{\longrightarrow}} \max  _ {\theta} \psi (\theta,\tilde{\tilde{\el}}(\theta),\tilde{\tilde{\eb}}(\theta)),
\end{equation}
with $\tilde{\tilde{\el}}(\theta)$ and $\tilde{\tilde{\eb}}(\theta)$ the solutions of score equations (\ref{eqz1}). By similar calculations as in the proof of Theorem \ref{theo2}, taking into account the Step  1, we can show that, for $\theta= \lim  _{n\rightarrow \infty} \theta_{nk}$, we have 
\begin{eqnarray*}
\Lambda_{nk} (\theta_{nk}) &=& \frac {1}{n} \sum  _{i \in I} \log \big[\theta +
\tilde{\tilde{\el}}^t(\theta) \eg_i(\tilde{\tilde{\eb}}(\theta))\big] + \frac {1}{n} \sum  _{j
\in J} \log  \big[(1-\theta) - \tilde{\tilde{\el}}^t (\theta)\eg_j(\tilde{\tilde{\eb}}(\theta))\big] \\ &&- \theta \log
\theta - (1-\theta) \log(1-\theta) + o_{\eP}(1).
\end{eqnarray*}
The above equation, together  with the law of large numbers, imply that $ \Lambda_{nk}(\theta_{nk}) \overset{{a.s.}} {\underset{n \rightarrow \infty}{\longrightarrow}} \psi( \theta, \tilde{\tilde{\el}}(\theta),
\tilde{\tilde{\eb}}(\theta)) $, where $\theta= \lim  _{n\rightarrow \infty} \theta_{nk}$. \ \\
For $\theta  \notin \{0,1,\theta_0\} $, partial derivative $ \partial \psi(\theta,\el,\eb)/\partial \theta$ becomes  \\
 $\displaystyle \int  _ {\Upsilon}  \int  _ {\eR}	\Big[ [\log K(\xx,e,\theta) \e1_{\{\theta <\theta_0\}}
dF_{\xx}(e)  + \log K(\xx,e,\theta) \e1_{\{\theta > \theta_0\}} dG_{\xx}(e) ] \ \\- [ \log (1-K(\xx,e,\theta)) \e1_{\{\theta
<\theta_0\}} dF_{\xx}(e) + \log (1-K(\xx,e,\theta)) \e1_{\{\theta > \theta_0\}} dG_{\xx}(e)] \Big]dH(\xx) \ \\  +\log (1-\theta) -\log \theta.$ \ \\ 
On the other hand, we have that, $dR_{\xx}{(e)}= \e1_{\{\theta < \theta_0\}} dF_{\xx}(e) +\e1_{\{\theta > \theta_0\}} dG_{\xx}(e)$. Hence, 
\begin{eqnarray*}
\frac{\partial \psi(\theta,\el,\eb)}{\partial \theta}&=& \int  _ {\Upsilon} \int  _ {\eR} [\log K(\xx,e,\theta) - \log (1-K(\xx,e,\theta))] dR_{\xx}{(e)} dH(\xx) \\&&
 -\log \theta +\log(1-\theta).
\end{eqnarray*}
 Because $\tilde{\tilde{\el}}^t(\theta) \eg(\xx,\tilde{\tilde{\eb}}(\theta))=K(\xx,e,\theta)-\theta$ and $ \ez_1(\theta,\tilde{\tilde{\el}}(\theta),\tilde{\tilde{\eb}}(\theta))=\textbf{0}_d$, we obtain 
\begin{equation}
\label{eqa}
\int  _ {\Upsilon} \big[ \int  _ {\eR} (1-\frac{\theta}{K(\xx,e,\theta)})dP_{\xx}{(e)}+\int  _ {\eR} (1-\frac{1-\theta}{1-K(\xx,e,\theta)})dQ_{\xx}{(e)} \big] dH(\xx)=0.
\end{equation}
On the other hand, we have $\ez_2(\theta,\tilde{\tilde{\el}}(\theta),\tilde{\tilde{\eb}}(\theta))=\textbf{0}_d$. Then
\begin{equation*}
\label{eqc}
\int  _ {\Upsilon}  \big[ \egg(\xx,\tilde{\tilde{\eb}}(\theta))  \int  _ {\eR} \frac{dP_{\xx}{(e)}}{K(\xx,e,\theta)}\big] dH(\xx)= \int  _ {\Upsilon}  \big[ \egg(\xx,\tilde{\tilde{\eb}}(\theta)) \int  _ {\eR} \frac{dQ_{\xx}{(e)}}{1-K(\xx,e,\theta)}\big] dH(\xx)=\textbf{0}_d.
\end{equation*}
Since  $\int  _ {\Upsilon}  \big[ \int  _ {\eR} dP_{\xx}(e) +\int  _ {\eR} dQ_{\xx}{(e)} \big ] dH(\xx)=1$, relation (\ref{eqa}) becomes 
\begin{equation}
\label{eqd}
1-\theta \int  _ {\Upsilon} \int  _ {\eR} \big[  \frac{dP_{\xx}(e)}{K(\xx,e,\theta)}- \frac{dQ_{\xx}{(e)}}{1-K(\xx,e,\theta)}\big] dH(\xx)- \int  _ {\Upsilon} \int  _ {\eR} \frac {dQ_{\xx}{(e)}}{1-K(\xx,e,\theta)}dH(\xx)=0.
\end{equation}
This relation is true for all $\theta \in (0,1)$. If we take $\theta=0$ and afterward $\theta=1$, relation  (\ref{eqd}) implies 
\begin{equation}
\label{eqe}
\int  _ {\Upsilon} \int  _ {\eR} \frac{dP_{\xx}{(e)}}{K(\xx,e,\theta)} dH(\xx)=\int  _ {\Upsilon} \int  _ {\eR} \frac{dQ_{\xx}{(e)}}{1-K(\xx,e,\theta)} dH(\xx)=1.
\end{equation}
The relation (\ref{argmax2}) is proved in a similar way as the proof of Theorem 3.2 of \cite{Liu:Zou:Zhang:08}, using relations (\ref{eqd}) and (\ref{argmax3}).\ \\
 {\it Step 3.} Similar as in the proof of Theorem 3.2 of \cite{Liu:Zou:Zhang:08}, we prove that $\psi (\theta_0,\tilde{\tilde{\el}}(\theta_0),\tilde{\tilde{\eb}}(\theta_0))$ \\$=0$, and that for all $\gamma \in (0, \min (\theta_0, 1-\theta_0)$, we have 
\begin{equation}
\label{argmax3}
\max   _{|k-n\theta_0| \geq n \gamma} \Lambda_{nk}(\theta_{nk}) \overset{{a.s.}} {\underset{n \rightarrow \infty}{\longrightarrow}}\max  _ {|\theta-\theta_0|\geq \gamma} \psi (\theta,\tilde{\tilde{\el}}(\theta),\tilde{\tilde{\eb}}(\theta)).
\end{equation}
 Which implies $\lim  _{n\rightarrow \infty} \eP [ | \argmax   _{k} \Lambda_{nk}(\theta_{nk}) - \theta_0 | \geq \gamma]=0$. \hspace*{\fill}$\blacksquare$  \\

\begin{Corollary}
\label{corol2}
The proof of Theorem \ref{theo4} implies that the estimator of $\theta_0$ defined by $$\tilde{\theta}_n\equiv  \min \{ \tilde{\theta}_{nk}; \tilde{\theta}_{nk}= \argmax   _{\theta_{nk} \in [\Theta_{1n}, \Theta_{2n}]} T_{nk}(\theta_{nk}, \tilde \el(\theta_{nk}), \tilde \eb(\theta_{nk})) \}$$
satisfies the property that $\tilde{\theta}_n - \theta_{n0} \overset{{\eP}} {\underset{n \rightarrow \infty}{\longrightarrow}} 0$. Taking into account Remark \ref{rem2}, we have also $\hat \theta_n -\theta_{n0} \rightarrow 0$ in probability, with $\hat \theta_n$ the estimator of $\theta_{0}$ defined in Remark \ref{rem2}. 
\end{Corollary}

 We prove by the following theorem that test statistic $\tilde{T}_n$ given by (\ref{eq8}) has the asymptotic power equal to 1.

%********************************************************************
%********************************************************************%********************************************************************%********************************************************************%********************************************************************%********************************************************************%********************************************************************

\begin{theorem}
\label{theo5}
 \noindent Under assumptions (A1)-(A7), the power of the empirical likelihood ratio  test $\tilde{T}_n$ converges to 1.
\end{theorem}

\noindent {\it Proof.} By the proof of Theorem \ref{theo4}, we have under hypothesis $H_1$
\begin{equation*}
\max  _ {\theta_{nk} \in [\Theta_{1n},\Theta_{2n}] }  \frac{T_{nk}(\theta_{nk}, \tilde \el(\theta_{nk}), \tilde \eb(\theta_{nk}))}{2n} \overset{{a.s.}} {\underset{n \rightarrow \infty}{\longrightarrow}} 0.
\end{equation*}
Taking into account relation (\ref{lbtild}), we have 
\begin{equation}
\label{rela}
\frac{\tilde{T}_n}{2n} \overset{{a.s.}} {\underset{n \rightarrow \infty}{\longrightarrow}} 0.
\end{equation}
In the other hand, if we suppose that the hypothesis $H_0$ is true, by Theorem \ref{theo3}, we have for all $ t \in \eR$,
\begin{equation*}
\eP \big [ A(\log u(n))\tilde{T}_n^{\frac{1}{2}} \leq t+ D(\log u(n))\big ] = exp(-e^{-t}).
\end{equation*}
Taking in the last relation $t=- \log \log u(n)$, we obtain 
\begin{equation}
 \lim  _{n\rightarrow \infty}  \eP \big [ A(\log u(n)) (\frac{\tilde{T}_n}{n})^{\frac{1}{2}} \leq -\frac{\log \log u(n)}{\sqrt{n}}+ \frac{D(\log u(n))}{\sqrt{n}}\big ] = 0.
\end{equation}
 The theorem follows. \hspace*{\fill}$\square$ \\

We emphasise that, similar results to Theorems \ref{theo3}, \ref{theo4} and \ref{theo5} were obtained by other authors for simpler models. The reader can find the corresponding results for the test to detect a change in distribution sequence in
\cite{Zou:Liu:Qin:Wang:14} and for detecting a change in the parameters of a linear model in \cite{Liu:Zou:Zhang:08}.
%********************************************************************
%********************************************************************%********************************************************************%********************************************************************%********************************************************************%********************************************************************%********************************************************************

\section{Extension to a particular two change-points model}
\label{sec2:1} 
\hh In this section, we consider the epidemic model. The epidemic linear model by empirical likelihood test was considered in paper \cite{Ni:pai:gup}. In a previous paper, \cite{yao} detect an epidemic alternative in the mean value of a sequence of independent normal random variables by various test statistics : likelihood ratio, recursive residual, score-like, semi-likelihood ratio. The works \cite{Ram:gup} and \cite{Ram} studied by likelihood ratio, the epidemic changes in the mean of a sequence of exponential random variables and respectively, in the shape parameter of a sequence of gamma random variables. \\

We assume under alternative hypothesis, denoted $H_2$, that the model have two change-points $k_1$ and $k_2$ ($1<k_1<k_2<n$), such that the model of the first and the third segment is the same. More specifically, the regression model can be written 
\begin{eqnarray}
\label{H2}
H_2 :      \textbf{Y}_i =
 \left\{
          \begin{array}{ll}

f(\textbf{X}_i, \eb_1)+ \varepsilon _i \,\,\,\, i=1,\cdots, k_1\\
           f(\textbf{X}_i,
\eb_2)+ \varepsilon _i \,\,\,\, i=k_1+1,\cdots,k_2 \\
f(\textbf{X}_i, \eb_1)+ \varepsilon _i \,\,\,\, i=k_2+1,\cdots, n.
          \end{array}
          \right.
\end{eqnarray} 
Therefore, we want to test the null hypothesis $H_0$ of no-change given by (\ref{eqq1}), against the alternative hypothesis $H_2$ given by (\ref{H2}).\ \\
\hh Under the hypothesis $H_2$, we consider the following two sets, $I'=\{1, ..., k_1, k_2+1, ..., n \}$ and $J'=\{ k_1+1, ..., k_2 \}$. We define the corresponding probability vectors $(u_1\cdots,u_{k_1},u_{k_2+1},\cdots, u_n)$ and
$(v_{k_1+1},\cdots,v_{k_2})$, where $u_i\equiv P [Y_i=y_i]$ and $v_j\equiv P[Y_j=y_j]$ denotes the probability to observe the value $y_i$ (respectively $y_j$), for the dependent variable $Y_i$ (respectively $Y_j$), for $i \in I'$ and $j \in J'$. Obviously, these probabilities satisfy the relations $\sum  _{i \in I'}u_i = 1$ and $\sum  _{j \in J'}v_j = 1$.

 Under  hypothesis $H_0$, the profile EL ratio for $\eb$ is  
\begin{eqnarray*}
 \textit{U}'_{0,n,k_1,k_2}(\eb)&=&
 \sup  _{(u_1,\cdots,u_{k_1},
u_{k_2+1},\cdots, u_n)}  \sup
 _{(v_{k_1+1},\cdots,v_{k_2})} \Big\{
 \prod  _{i \in I'} (n-k_2+k_1)u_i
 \prod  _{j \in J' }(k_2-k_1)v_j ; \\ &&
 \sum   _{i\in I'} u_i =
\sum   _{j\in J'} v_j =1, \sum   _{i\in
I'} u_i \eg_i(\eb )=
 \sum  
_{j\in J'} v_j  \eg_j( \eb )=\textbf{0}_d \Big\}.
\end{eqnarray*}

Under hypothesis $H_2$, the profile EL ratio for $\eb_1$, $\eb_2$ has the form 
\begin{eqnarray*}
\textit{U}'_{1,n,k_1,k_2}(\eb_1, \eb_2)&= &
\sup  _{(u_1,\cdots,u_{k_1}, u_{k_2+1},\cdots, u_n)}
 \sup  _{(v_{k_1+1},\cdots,v_{k_2})} 
\Big \{  \prod  _{i\in I'}(n-k_2+k_1)u_i
\prod  _{j\in J'}(k_2-k_1)v_j  
;\\ && \sum   _{i\in I'} u_i
\textbf{g}_i (\eb_1 )=  \sum
 _{j\in J'} v_j\textbf{g}_j (\eb_2 )=\textbf{0}_d  \Big\}.
 \end{eqnarray*}
Then, in order to test $H_0$ against $H_2$, we consider  the profile EL ratio  $
\textit{U}'_{0,n,k_1,k_2}(\eb)/\textit{U}'_{1,n,k_1,k_2}(\eb_1,\eb_2) $.\\

Similarly as in Section \ref{sec1:1},  when we tested a single change-point, using Lagrange multipliers, we obtain that under hypothesis $H_0$, the probabilities $u_i$, $v_j$ are
\begin{equation}
\label{p4}
 u_i= \frac{1}{ (n-k_2+k_1)+  n \el^t_1\eg_i(\eb )} \qquad  v_j= \frac{1}{(k_2-k_1) - n \el^t_2
\textbf{g}_j(\eb )}.
\end{equation}

 Using the similar arguments as in the proof of Theorem \ref{theo1}, we deduce that the asymptotic distribution of $- 2 \log U'_{1,n,k_1,k_2}(\eb_1,\eb_2)$ is $\chi^2(3d)$ and then we can  consider the test statistic $- 2 \log U'_{0,n,k_1,k_2}(\eb)$. 
 We restricted to the case where $\el_1$ and $\el_2$ satisfy the constraint $\tilde{\eVV}_{1n}(\eb) \el_1= \tilde{\eVV}_{2n}(\eb) \el_2$, with  $\tilde{\eVV}_{1n}(\eb)=(n+k_1-k_2)^{-1} \sum   _{i\in I'} \egg(\eb)$ and  $\tilde{\eVV}_{2n}(\eb)=(k_2-k_1)^{-1} \sum   _{j\in J'} \egg(\eb)$.  In this case, considering the parameter $\theta_{n,k_1,k_2}=n^{-1}(n-k_2+k_1)$, that depends on two change-points $k_1$, $k_2$, we will consider the test statistic 
\begin{equation}
\label{p8}
2 \sum
  _{i \in I'}\log \Big[{1+\frac{1}{\theta_{n,k_1,k_2}}
\el^t\eg_i( \eb )}\Big]+ 2\sum   _{j \in
J'}\log \Big[{1- \frac{1}{1-\theta_{n,k_1,k_2}}\el^t \tilde{\eVV}_{1n}(\eb)
\tilde{\eVV}_{2n}^{-1}(\eb)\eg_j (\eb )}\Big].
\end{equation}
Let us denote by $\hat{\el}(\theta_{n,k_1,k_2})$, $\hat{\eb }(\theta_{n,k_1,k_2})$ the solutions of the score equations of this random process equal to zero.  We can show, as in Section 2, that  statistic (\ref{p8}) is, under hypothesis $H_0$, asymptotically equivalent to the statistic 
\begin{equation*}
U_{n,k_1,k_2}(\theta_{n,k_1,k_2},\el,\eb) \equiv 2 \sum
  _{i \in I'}\log \Big[{1+\theta_{n,k_1,k_2}^{-1}
\el^t\eg_i( \eb )}\Big]+ 2\sum   _{j \in
J'}\log \Big[{1-(1-\theta_{n,k_1,k_2})^{-1}\el^t \eg_j (\eb )}\Big].
\nonumber
\end{equation*}
Then, we will consider for testing  null hypothesis $H_0$ against $H_2$ the test statistic 
 $$\max  _{1< k_1<k_2<n}  \{ U_{n,k_1,k_2}(\theta_{n,k_1,k_2},\hat{\el}(\theta_{n,k_1,k_2}),\hat{\eb }(\theta_{n,k_1,k_2})) \}.$$
\hh  In the case when $k_1$ or $k_2-k_1$ have a small value, the maximum  empirical likelihood may not exist. In this case, the proposed test may not detect the presence of change in the
model. For the empirical likelihood maximum always exists, we consider two natural numbers $\Theta_{n1}$ and $\Theta_{n2}$, such as $\Theta_{n1}<k_1<k_2<n-\Theta_{n2}$. Finally, the test statistic for testing $H_0$ against $H_2$ becomes 
$$\max  _{\Theta_{n1}<k_1<k_2<n-\Theta_{n2}}  \{ U_{n,k_1,k_2}(\theta_{n,k_1,k_2},\hat{\el
}(\theta_{n,k_1,k_2}),\hat{\eb }(\theta_{n,k_1,k_2})) \}.$$
\hh  In order to facilitate the practical utilization of the test statistic, we can easily obtain the corresponding statistic given in Remark \ref{rem2} by relation (\ref{eqdez}).

\section{Simulation study}
\label{sec3:1}
\hh In this section, we report a simulation study by Monte Carlo method, in order to evaluate the performance of the proposed test statistics. Firstly, for a fixed theoretical size, we calculate the critical value of the test statistics, for different values of $n$. Afterward, by Monte Carlo technique, we calculate empirical test size, empirical test power and the estimation of the change-point localisation. This study was conducted firstly for a nonlinear model with a single change-point and secondly for an epidemic model. The obtained results by the proposed test statistic are compared with whose obtained from LS method, proposed by \cite{Boldea:Hall:13}.\\
\hh All simulations were performed using the R language. The program codes are available from the authors. \\
We consider the nonlinear function 
\begin{equation}
\label{model1}
f(x,\eb)=a\frac{1-x^b}{b},
\end{equation}
 with $\eb=(a,b) \in [-100, 100] \times [0.1, 20]$. The same model was considered in  \cite{Ciuperca:11b}, where the model was estimated by the penalized least absolute deviation method. 
\subsection{Model with a single change-point}
\hh  For the nonlinear function of (\ref{model1}), the following two-phase (one change-point) nonlinear model is considered under $H_1$
\begin{equation}
\label{model2}
Y_i =  a_1 \frac{1-X_i^{b_1}}{b_1} \e1_{i \leq k_0}+ a_2 \frac{1-X_i^{b_2}}{b_2} \e1_{i >k_0}+ \varepsilon_i, \qquad i=1, \cdots , n
\end{equation}
\noindent with $X_i=i/1000$, $n=1000$ and  true value of parameters $a_1^0=10$, $b_1^0=2$, $a^0_2=7$, $b_2^0=1.75$. Under hypothesis $H_0$, the true parameters are $a^0=10$, $b^0=2$. \\
The change absence against one-change in model is tested using the (approached) maximum empirical likelihood statistic $T(\theta_{nk})$ given by (\ref{eqdez}). \ \\
\hh  In order to calculate the empirical test size, an without change-point model is considered and we count, the number of times, on the Monte Carlo replications when we obtain $ \max   _{\theta_{nk}} T^ {1/2}(\theta_{nk}) \geq c_{\alpha}$. For a fixed size $\alpha \in (0,1)$, critical value $c_{\alpha}$ is calculated in accordance with Corollary \ref{corol1} : $$c_{\alpha}= \frac{-\log (-\log \alpha)+D(\log u(n))}{A(\log u(n))}.$$
For theoretical size $\alpha=0.05$, we first calculate critical values $c_\alpha$, varying the sample size $n$ from 200 to 1000 (see Table 1).\\
\hh For model (\ref{model2}) with Gaussian standardized errors, 500 Monte Carlo replications were performed.  We also present in Table 1 the empirical power, using statistic test (\ref{eqdez}). For different position of change-point. For any change-point location, the asymptotic test power is 1.   
\begin{table}[h]
\begin{center}
\caption{Critical values $c_\alpha$, for   $\alpha$=0.05. Empirical power on 500 Monte Carlo replications, when $\varepsilon \sim {\cal N}(0,1)$.}
\begin{tabular}{|cccc|} \hline \noalign{\smallskip}
n & $k_0$ &  $c_{\alpha}$  & power \\  
\noalign{\smallskip}\hline\noalign{\smallskip}
1000 & 600 &1.544 & 1  \\
  800 & 500 & 1.492  & 1 \\
   600 & 400 & 1.434  & 1 \\
 400 & 250 & 1.340  & 1  \\
 200 & 75 & 1.133  & 1  \\ 
 \noalign{\smallskip}\hline
\end{tabular}
\end{center}
\label{tabl1}
\end{table} 
We fix sample size $n=1000$, theoretical test size  $\alpha=0.05$ and we vary the error distribution. In order to calculate the empirical size of test (type I error probabilities),  500  Monte Carlo replications are realized for different error distributions: $\varepsilon_i={\cal N}(0,1)$, $\varepsilon_i=2{\cal E}xp(2)-1$, $\varepsilon_i= 1/\sqrt{6}(\chi^2(3)-3)$ and $\varepsilon_i=2/\sqrt{6}   t(6)$, where  ${\cal N}(0,1)$,  ${\cal E}xp(2)$,  $\chi^2(3)$ and $t(6)$ are standard normal distribution, exponential distribution with  mean 1/2, chi-square distribution with degree of freedom  3 and Student distribution with degree of freedom  6, respectively. In all cases, except for Student distribution (when the empirical size is slightly larger than 0.05), the empirical size is 0 (see Table \ref{tabl3:1}).
For the same four error distributions, but for model with a change-point in $k_0$, by  500 Monte Carlo model replications, for different change-point location: $k_0 \in \{200, 400, 600, 800 \}$, we obtain that the empirical power is 1, in any case. 
\begin{table}[h]
\begin{center}
\caption{Empirical size for four error distributions on 500 Monte Carlo replications, $\alpha$=0.05.}
\label{tabl3:1}
\begin{tabular}{| c | c  |c |  c| c|} \hline \noalign{\smallskip}
 $n$ &  \multicolumn{4}{c|}{ error distribution} \\
 &  Normal &Exponential &$\chi^2$  & Student \\ 
\noalign{\smallskip}\hline\noalign{\smallskip} 
 200 & 0& 0&0&0 \\ 
 1000 & 0& 0&0&0 \\
  \noalign{\smallskip}\hline
\end{tabular}
\end{center}
\end{table} 
\begin{table}[h]
\caption{Descriptive statistics for the estimators of the change-point. Model with two phases by EL method, $n=1000$, 500 Monte Carlo replications. }
%\begin{center}
\begin{tabular}{|c|c|c|c|c|c|c|}\hline \noalign{\smallskip}
  error distribution & $k^0$  & \multicolumn{5}{c|}{ $\hat k_n $}  \\ 
    \cline{3-7} 
     & & min($\hat k_n $) & max($\hat k_n $) & mean($\hat k_n $)  & sd($\hat k_n $) & median($\hat k_n $) \\ \noalign{\smallskip}\hline\noalign{\smallskip}
      $\varepsilon_i \sim {\cal N}(0,1)$ & 200 & 186 & 205 & 197 & 4 & 199 \\
      & 400 & 388 & 416 & 400 & 4 & 400 \\
      & 600 & 585 & 612 & 600 & 6 & 600 \\
      & 800 & 780 & 819 & 795 & 9 & 797 \\ 
            \noalign{\smallskip}\hline
 $\varepsilon_i \sim 2/\sqrt{6} t(6)$ &  200 & 180 & 209 & 199 & 5 & 200 \\
      & 400 & 390 & 407 & 400 & 4 & 400 \\
      & 600 & 580 & 611 & 599 & 6 & 600 \\
      & 800 & 780 & 815 & 797 & 10 & 798 \\
            \noalign{\smallskip}\hline
 $\varepsilon_i \sim 2{\cal E}xp(2)-1$&  200 & 183 & 204 & 199 & 4 & 200 \\
      & 400 & 391 & 412 & 401 & 5 & 400 \\
      & 600 & 591 & 609 & 600 & 4 & 600 \\
      & 800 & 780 & 817 & 795 & 5 & 799 \\
            \noalign{\smallskip}\hline
 $\varepsilon_i \sim 1/\sqrt{6}(\chi^2(3)-3)$&  200 & 188 & 210 & 200 & 4 & 200 \\
      & 400 & 393 & 417 & 401 & 5 & 400 \\
      & 600 & 585 & 614 & 600 & 6 & 600 \\
      & 800 & 780 & 807 & 795 & 8 & 797 \\ 
       \noalign{\smallskip}\hline
\end{tabular}
%\end{center}
\label{Tabl6:1} 
\end{table} \\
\hh As mentioned in Remark \ref{rem2}, one can also estimate the change-point location by EL method. In table \ref{Tabl6:1} we have the summarized results (minimum, maximum, mean, standard-deviation, median) for the estimator $\hat k_n $, given by Remark \ref{rem2}, by 500 Monte Carlo replication. In view of the results presented in Table \ref{Tabl6:1}, for different error distribution and for  different positions of the change in the interval, we deduce that the proposed estimation method approaches very well the true value $k^0$, regardless of the error distribution and of the change-point position on the interval $[1:n]$.  Note that, in all situations the median and the mean of the change-point estimations coincide or is very close to the true value.
\subsection{Epidemic model}
\hh For nonlinear function of (\ref{model1}), under hypothesis $H_2$, we consider the following three-phase (two change-points) model 
\begin{equation}
\label{model3}
Y_i = a_1 \frac{1-X_i^{b_1}}{b_1}\e1_{i \leq k_1}+a_2 \frac{1-X_i^{b_2}}{b_2}\e1_{k_1<i \leq k_2}+   a_1 \frac{1-X_i^{b_1}}{b_1} \e1_{k_2 <i \leq n}+\varepsilon_i,
\end{equation}
\noindent with $X_i=i/1000$, $n=1500$ and the true value of  parameters $a_1^0=10$, $b_1^0=2$, $a^0_2=7$, $b_2^0=1.75$. Under null hypothesis $H_0$ the true parameters are $a^0=10$, $b^0=2$. \\
\hh In Table \ref{Tabl5:1} we give results after 150 Monte Carlo replications in order to calculate the empirical power of test, for $n=1500$. We deduce that empirical size is zero and empirical test power is 1.\\
\begin{table}[h]
\caption{Empirical powers and empirical size for epidemic model, $\alpha=0.05$, $n=1500$.}
\begin{center}
\label{Tabl5:1}
\begin{tabular}{|c c c |} \hline \noalign{\smallskip}
$k_1$ &$k_2$ & power  \\  
\noalign{\smallskip}\hline\noalign{\smallskip}
\multicolumn{2}{|c}{no-change} &0 \\
 100 & 900  & 1 \\
200 & 500  & 1   \\
400 & 600  & 1  \\
600 & 900  & 1  \\ 
 \noalign{\smallskip}\hline
\end{tabular}
\end{center}
\end{table} 
 \subsection{Comparison with LS test}
\hh On data considered in subsection 4.1 for $\varepsilon \sim {\cal N}(0,1)$ and $n=1000$ we apply the method proposed by \cite{Boldea:Hall:13}, where the estimation method and the associated test is by least squares. This study is realized by computing the test statistic  $\sup F(0:1)$ given in \cite{Boldea:Hall:13}. Under hypothesis $H_1$, when the model has a change-point in $k^0=600$, we realize 500  Monte Carlo  simulations. We obtain that the test statistic value always exceeds the critical value of  12.85 (see \cite{Bai:Perron}). Then, by LS test of \cite{Boldea:Hall:13}, the null hypothesis  $H_0$ is always rejected and hence the power of test is 1. Whereas if we generate the values $Y_i$ without change-point for gaussian errors, then, the test statistic value of \cite{Boldea:Hall:13} always exceeds critical. Hence the empirical size of the test proposed by \cite{Boldea:Hall:13} is 1, a result significantly worse than that obtained by our test. 
We note that (see Table \ref{Tabl7:1}) if under $H_1$ the true change-point is off-centered in the measurement interval, because of the function nonlinearity, then numerical problem arise for the LS estimation method. This is symbolized by "???" in Table \ref{Tabl7:1}. The same problem appears when the errors are not gaussian, regardless of the position of the change-point in the measurement interval. In contrast, we have seen that the EL test works for any error distribution and any change-point position.
\begin{table}
\caption{Descriptive statistics for the estimators of the change-point. Model with two phases by LS method, $n=1000$, 500 Monte Carlo replications. }
%\begin{center}
\begin{tabular}{|c|c|c|c|c|c|c|}\hline \noalign{\smallskip}
  error distribution & $k^0$  & \multicolumn{5}{c|}{ $\hat k_n $}  \\ 
    \cline{3-7} 
     & & min($\hat k_n $) & max($\hat k_n $) & mean($\hat k_n $)  & sd($\hat k_n $) & median($\hat k_n $) \\ \noalign{\smallskip}\hline\noalign{\smallskip}
      $\varepsilon_i \sim {\cal N}(0,1)$ & 200 & ??? & ??? & ???  & ???  & ???  \\
      & 400 & 396 & 400 & 399 & 1 & 400 \\
      & 600 & 595 & 605 & 600 & 3 & 600 \\
      & 800 & ??? & ??? & ??? & ??? & ??? \\ 
            \noalign{\smallskip}\hline
 \end{tabular}
%\end{center}
\label{Tabl7:1} 
\end{table}

\section{Appendix}
\label{sec4:1}

%++++++++++++++++++++++++++++++++++++++++++++++++++++++++++++++++++++
%++++++++++++++++++++++++++++++++++++++++++++++++++++++++++++++++++++
%++++++++++++++++++++++++++++++++++++++++++++++++++++++++++++++++++++
%++++++++++++++++++++++++++++++++++++++++++++++++++++++++++++++++++++
%++++++++++++++++++++++++++++++++++++++++++++++++++++++++++++++++++++

The following lemma will be used in the proof of propositions, theorems and of other lemmas.

\begin{lemma}
\label{lem0}
\noindent Let  $\XX=(X_1,\cdots,X_p)$ a random vector (column), with the random variables $X_1,\cdots,X_p$ not necessarily independent, and $\textbf{M}=(m_{ij})_{1 \leq i,j \leq p}$, such that $\textbf{M}= \XX \XX^t$. If for j=1, ..., p, we have
  \begin{equation}
\label{eqeta}
for\,\,all \,\,  \eta_j>0, there\,\,exists \,\, \delta_j >0\,\,\,\,\, such\,\, that \,\,\,\, \eP [|X_j|\geq \delta_j] \leq \eta_j, 
\end{equation}
then \ \\
 $(i)\,\, \eP \big[ \| \XX \|_1 \geq p \max  _{1 \leq j \leq p}\delta_j \big]\leq \max  _{1 \leq j \leq p}  \eta_j$, \\
 $(ii)\,\, \eP \big[ \| \XX \|_2 \geq \sqrt{p} \max  _{1 \leq j \leq p}\delta_j\big]\leq \max   _{1 \leq j \leq p}\eta_j$,  \\
 $(iii) \,\, \eP \big[ \| \textbf{M} \|_1 \geq p \max   _{1 \leq i,j \leq p} \{\delta_i^2,\delta_j^2\}\big]\leq \max   _{1 \leq i,j \leq p} \{\eta_i^2,\eta_j^2\}$, \ \\
 where $\| \textbf{M} \|_1= \max  _{1 \leq j \leq p} \{ \sum   _{i=1}^p |m_{ij}| \}$ is the subordinate norm to the vector norm $\|.\|_1$. 
 \end{lemma}
 
\noindent {\it Proof of Lemma \ref{lem0}}. $(i)$ Using the relation (\ref{eqeta}), we can write 
\begin{equation*}
  \eP [\| \XX \|_1 \geq  p \max  _{1 \leq j \leq p} \delta_j ] 
  \leq  \eP [p \max  _{1 \leq j \leq p}  | X_j | \geq p \max  _{1 \leq j \leq p}  \delta_j] \leq \max  _{1 \leq j \leq p}   \eta_j.
 \end{equation*}
 $(ii)$ The relation (\ref{eqeta}) is equivalent to $\eP \big[X_j^2\geq \delta_j^2\big] \leq \eta_j$, which implies that 
\begin{equation*}
  \eP [\| \XX \|_2^2 \geq  p \max  _{1 \leq j \leq p}  \delta_j^2] 
= \eP [\max  _{1 \leq j \leq p} X_j^2 \geq  \max  _{1 \leq j \leq p} \delta_j^2] \\ 
  \leq \max  _{1 \leq j \leq p}  \eta_j.
 \end{equation*}
  $(iii)$ For $1\leq i,j \leq p$, we have 
\begin{eqnarray*}
  \eP [| X_i X_j| \geq  \max    \{ \delta_i^2,\delta_j^2\} ] 
  \leq \eP [\max  \{X_i^2,X_j^2\}\geq  \max  \{\delta_i^2,\delta_j^2 \} ] 
\leq \max   \{\eta_i^2,\eta_j^2\}.
 \end{eqnarray*}
  Then, $\eP [| m_{ij} | \geq  \max  \{\delta_i^2,\delta_j^2\}] 
\leq \max   \{\eta_i^2,\eta_j^2\}$. Hence, for each $1\leq j \leq p$, 
\[
\eP [\sum  _ {i=1} ^p| m_{ij} | \geq p \max   \{\delta_i^2,\delta_j^2\}\big] \leq \eP \big[p \max   _{1 \leq i \leq p}|m_{ij} | \geq p \max   \{ \delta_i^2,\delta_j^2\}]
\leq \max \{\eta_i^2,\eta_j^2\}.
\]
\hspace*{\fill}$\blacksquare$ \\
%%%%%%%%%%%%%%%%%%%%%%%%%%%%%%%%%%%%%%%%%%%%%%%%%%%%%%%%%%%%%%%%%%%%%%%%%%%%
%%%%%%%%%%%%%%%%%%%%%%%%%%%%%%%%%%%%%%%%%%%%%%%%%%%%%%%%%%%%%%%%%%%%%%%%%%%%
%%%%%%%%%%%%%%%%%%%%%%%%%%%%%%%%%%%%%%%%%%%%%%%%%%%%%%%%%%%%%%%%%%%%%%%%%%%%

\begin{lemma}
\label{lem3}
Let the $\eta$-neighbourhood of $\ebo$, ${\cal V}_{\eta}(\ebo)= \{ \eb \in \Gamma; \| \eb-\ebo\|_2 \leq \eta \}$, with $\eta \rightarrow 0$. Then, under assumptions (A1)-(A4), for all $\epsilon>0$, there exists a positive constant $M>0$, such that, for all $\eb \in {\cal V}_{\eta}(\ebo)$,
\begin{equation*}
\label{eqg}
\eP \big[\|\eg_i(\eb) \|_1 \geq M \big] \leq \epsilon.
\end{equation*}
\end{lemma}
 \noindent {\it Proof of Lemma \ref{lem3}}. In the following, for simplicity, we denote the functions $\ef(\XX_i,\eb)$ by $\ef_i(\eb)$, and $\eff(\XX_i,\eb)$ by $\eff_i(\eb)$.
The Taylor's expansion up the order 2 of $\eg_i(\eb)$ at $\eb=\eb^0$ is 
%\begin{eqnarray}
%\label{eq99}
%\eg_i(\eb)&=& \ef_i(\ebo) \varepsilon_i+\frac{1}{2}\eff_i(\tilde{\eb}_i) (\eb-\ebo) \varepsilon_i  -\frac{1}{2}\ef_i(\ebo)\ef_i^t (\ebo) (\eb-\ebo) \nonumber  \\ && 
% -\frac{1}{6}\ef_i(\ebo)(\eb-\ebo)^t \eff_i(\tilde{\tilde{\eb}}_i) (\eb-\ebo) - \frac{1}{4}\eff_i(\tilde{\eb}_i) (\eb-\ebo)\ef_i^t (\ebo) (\eb-\ebo) \nonumber\\ && 
%-\frac{1}{12}\eff_i(\tilde{\eb}_i) (\eb-\ebo)(\eb-\ebo)^t \eff_i(\tilde{\tilde{\eb}}_i)(\eb-\ebo) ,
%\end{eqnarray} 
\begin{eqnarray}
\label{eq99}
\eg_i(\eb)&=& \ef_i(\ebo) \varepsilon_i+\frac{1}{2}\eM_{1i} (\eb-\ebo) \varepsilon_i  -\frac{1}{2}\ef_i(\ebo)\ef_i^t (\ebo) (\eb-\ebo) \nonumber  \\ && 
 -\frac{1}{6}\ef_i(\ebo)(\eb-\ebo)^t \eM_{2i} (\eb-\ebo) - \frac{1}{4}\eM_{1i} (\eb-\ebo)\ef_i^t (\ebo) (\eb-\ebo) \nonumber\\ && 
-\frac{1}{12}\eM_{1i} (\eb-\ebo)(\eb-\ebo)^t \eM_{2i}(\eb-\ebo) ,
\end{eqnarray} 
where $\eM_{1i}=  \Bigg( \frac{\partial^2 f_i( \eb^{(1)}_{i,jk})}{\partial \beta_j \partial \beta _k} \Bigg)_{1\leq j,k \leq d} $, $\eM_{2i}=  \Bigg( \frac{\partial^2 f_i( \eb^{(2)}_{i,jk})}{\partial \beta_j \partial \beta_ k} \Bigg)_{1\leq j,k \leq d} $ and \\ $ \eb^{(1)}_{i,jk}=\ebo+u_{i,jk}(\eb- \ebo)$, $ \eb^{(2)}_{i,jk}=\ebo+v_{i,jk}(\eb- \ebo)$, with $u_{i,jk}, v_{i,jk} \in [0,1]$. %where, $\tilde{\eb}_i =\ebo+ \eu( \eb-\ebo)$, $\tilde{\tilde{\eb}}_i =\ebo+ \ev( \eb-\ebo)$, with $\eu,\ev \in [0,1]^d$. 
We note that $ \eb^{(1)}_{i,jk}$ and $ \eb^{(2)}_{i,jk}$ are random vectors which depend on $\XX_i$.\ \\
For $\ef_i(\ebo) \varepsilon_i$, because $\XX_i $ and $\varepsilon_i$ are independent, and $\eE(\varepsilon_i)=0$, we have that $\eE[\ef_i(\ebo) \varepsilon_i]=0$ and $\Var[\ef_i(\ebo) \varepsilon_i]= \sigma^2 \eV$.
 For the j-th component of $\ef_i(\ebo)$, by the Bienaymé-Tchebychev's inequality, for $1 \leq j \leq d$, for all $\epsilon_1>0$, we have 
\begin{equation}
\label{eqbt1}
\eP\big[ |\frac{\partial f_i(\ebo)}{\partial \beta_j} \varepsilon_i |\geq \epsilon_1 \big]\leq \frac{\sigma^2}{\epsilon_1^2} V_{jj},
\end{equation}
\noindent where $V_{jj}$ is the j-th term diagonal of the matrix $\eV$. \ \\
\noindent For all $ \epsilon>0$, taking $\epsilon_1= \sigma \sqrt{6V_{jj}/\epsilon}$ in (\ref{eqbt1}), we obtain $\eP \big  [ | \frac{\partial f_i(\ebo)}{\partial \beta_j}  \varepsilon_i | \geq \sigma \sqrt{6V_{jj}/\epsilon} \big] \leq \epsilon / 6$. Applying Lemma \ref{lem0} (i), we obtain, for all $ \epsilon>0$
\begin{equation}
\label{eqbt3}
\eP \big[ \| \ef_{i}(\ebo) \varepsilon_i \|_1 \geq {\frac{\sigma d}{\sqrt{\epsilon}}} \max   _{1\leq j \leq d} \sqrt{6V_{jj}} \big]\leq \epsilon / 6.
\end{equation}

%***********************************************************************************
%***********************************************************************************
For the second term of the right-hand side of (\ref{eq99}), using assumption (A3), we obtain that for $  1\leq j,k \leq d$, for all $\epsilon>0$ there exists $\epsilon_2>0$, such that, $
\eP \big[|\frac{\partial^2 f_i( \eb^{(1)}_{i,jk})}{\partial \beta_j \partial \beta_k} | \geq  \epsilon_2\big] \leq \epsilon/6$.
By Lemma \ref{lem0} (iii), we have that for all $\epsilon>0$, 
\begin{equation}
\label{eqbt14}
\eP \big[\| \eM_{1i}  \|_1 \geq  \epsilon_2 \big] \leq \frac{\epsilon}{6}.
\end{equation}
 Using Bienaymé-Tchebychev's inequality, and assumption (A1), we obtain that for all $C_1>0$
\begin{equation}
\label{epselon}
\eP \big[| \varepsilon_i | > C_1 \big] \leq  \frac{\sigma^2}{C_1}.
\end{equation}
Recall that $\| \eb-\ebo\|_2<\eta$, with $\eta \rightarrow 0$. Then, using (\ref{eqbt14}) and (\ref{epselon}), we can write that, for all $\epsilon>0$, there exists $\epsilon_2> 0$ such that, $
\eP \big[ \|\eM_{1i}  (\eb-\ebo) \varepsilon_i \|_1 \geq \epsilon_2 \big]
 \leq  \eP \big[ \|\eM_{1i}  \| _1 |\varepsilon_i| \,\, \|\eb-\ebo  \|_1 \geq \epsilon_2 \big] \leq \eP \big[ \| \eM_{1i}   \| _1 \geq \epsilon_2/C_1 \eta \big] 
\leq \eP \big[ \|\eM_{1i}  \| _1  \geq \epsilon_2 \big] 
\leq  \epsilon/6$.
 Therefore, for all $\epsilon>0$, there exists $\epsilon_2>0$ such that
  \begin{equation}
\label{eqbt15}
\eP \big[ \| \eM_{1i} (\eb-\ebo) \varepsilon_i \|_1 \geq \epsilon_2 \big]\leq  \frac{\epsilon}{6}.
\end{equation}

We consider now the term $\ef_i(\ebo)\ef_i^t (\ebo) (\eb-\ebo)$ of relation (\ref{eq99}). By Markov's inequality, taking also into account assumption (A4), we obtain for $1 \leq j,l \leq d$, for all $\epsilon_3>0$, that $\eP \big[ | \frac{\partial f_i(\ebo)}{\partial \beta_j}  \frac{\partial f_i(\ebo)}{\partial \beta_l}  | \geq \epsilon_3 \big] \leq \eE[| \frac{\partial f_i(\ebo)}{\partial \beta_j}  \frac{\partial f_i(\ebo)}{\partial \beta_l}  |]/\epsilon_3$.
 We choose, for all $\epsilon>0, \epsilon_3= 6 \eE[| \frac{\partial f_i(\ebo)}{\partial \beta_j} \frac{\partial f_i(\ebo)}{\partial \beta_l}  | ]/\epsilon$. Then, the last relation becomes $\eP \big[ | \frac{\partial f_i(\ebo)}{\partial \beta_j}  \frac{\partial f_i(\ebo)}{\partial \beta_l}  |  \geq  6 \eE[| \frac{\partial f_i(\ebo)}{\partial \beta_j}  \frac{\partial f_i(\ebo)}{\partial \beta_l}  |  \big] / \epsilon \big ] \leq  \epsilon/6$.
 Using Lemma \ref{lem0} (iii), we obtain
\begin{equation*}
\label{eqbt6}
\eP \bigg[ \| \ef_{i} (\ebo) \ef^t_{i} (\ebo) \|_1 \geq  \frac{6d}{\epsilon} \max   _{1\leq j,l \leq d}\eE[| \frac{\partial f_i(\ebo)}{\partial \beta_j}  \frac{\partial f_i(\ebo)}{\partial \beta_l}  | ] \bigg]\leq \frac{\epsilon}{6},
\end{equation*}
relation that involves, since for all $C_2 >0$ we have $\| \eb-\eb_0 \|_1 \leq C_2 \eta$ for $\eta \rightarrow 0 $, that \ \\
 $\eP \big[\| \ef_{i} (\ebo) \ef^t_{i} (\ebo) (\eb-\ebo) \|_1 \geq   6d/\epsilon \max   _{1\leq j,l \leq d}\eE[ | \frac{\partial f_i(\ebo)}{\partial \beta_j}  \frac{\partial f_i(\ebo)}{\partial \beta_l}  | ] \big]  
\ \\
\leq \eP \big[ \| \ef_{i} (\ebo) \ef^t_{i} (\ebo) \|_1 \geq  6d/\epsilon \max   _{1\leq j,l \leq d}\eE[ | \frac{\partial f_i(\ebo)}{\partial \beta_j}  \frac{\partial f_i(\ebo)}{\partial \beta_l}  | ] \big]  
\leq \epsilon/6. $ \ \\ 
 Then, for all $\epsilon>0$
\begin{equation}
\label{eqbt7}
\eP \big[\| \ef_{i} (\ebo) \ef^t_{i} (\ebo) (\eb-\ebo) \|_1 \geq   \frac{6d}{  \epsilon} \max   _{1\leq j,l \leq d} \eE[ | \frac{\partial f_i(\ebo)}{\partial \beta_j}  \frac{\partial f_i(\ebo)}{\partial \beta_l}  | ] \big] \leq \frac{\epsilon}{6}.
\end{equation}
\hspace*{0.5cm} For $\eM_{1i} (\eb-\ebo) \ef^t_{i} (\ebo) (\eb-\ebo)$ of relation (\ref{eq99}), using assumption (A3) and the Markov's inequality, we obtain for each j-th component $\frac{\partial f_i(\ebo)}{\partial \beta_j}$ of the vector $\ef_{i} (\ebo)$, for all $\epsilon_4>0$, that $
\eP \big[ |\frac{\partial f_i(\ebo)}{\partial \beta_j}| \geq \epsilon_4 \big] \leq \eE[|\frac{\partial f_i(\ebo)}{\partial \beta_j}|]/\epsilon_4$. We choose, for all $\epsilon>0$, $\epsilon_4=6\eE[|\frac{\partial f_i(\ebo)}{\partial \beta_j}|]/\epsilon$ and this last relation becomes $
\eP \big[ |\frac{\partial f_i(\ebo)}{\partial \beta_j}| \geq 6\eE[|\frac{\partial f_i(\ebo)}{\partial \beta_j}|] / \epsilon \big] \leq  \epsilon/6$. Applying Lemma \ref{lem0} (i), for all $\epsilon>0$ we obtain 
\begin{equation}
\label{eqbt12}
\eP \big[ \| \ef_{i} (\ebo) \|_1 \geq    \frac{ 6d}{\epsilon} \max   _{1\leq j \leq d} \eE[|\frac{\partial f_i(\ebo)}{\partial \beta_j}|]\big] \leq  \frac{\epsilon}{6}.
\end{equation}
 Using assumption (A3), and relations (\ref{eqbt14}), (\ref{eqbt12}),  we can write that  \ \\
 $\eP \big[\| \eM_{1i} (\eb-\ebo) \ef^t_{i} (\ebo) (\eb-\ebo) \|_1 \geq   6d/ \epsilon \max   _{1\leq j \leq d} \eE[|\frac{\partial f_i(\ebo)}{\partial \beta_j}|]\big] \ \\ \leq \eP \big[ \| \ef^t_{i} (\ebo)\|_1 \geq 6d/ \epsilon \max   _{1\leq j \leq d} \eE  [|\frac{\partial f_i(\ebo)}{\partial \beta_j}|]\big] 
\leq \epsilon/6 $.  \ \\ 
Therefore, for all $\epsilon>0$,
\begin{equation}
\label{eqbt13}
\eP \big[\| \eM_{1i} (\eb-\ebo) \ef^t_{i} (\ebo) (\eb-\ebo) \|_1 \geq \frac{ 6d}{\epsilon} \max   _{1\leq j \leq d} \eE[|\frac{\partial f_i(\ebo)}{\partial \beta_j}|]\big] \leq \frac{\epsilon}{6}.
\end{equation}
\hspace*{0.5cm} Taking into account assumptions (A3), (A4), by relations (\ref{eqbt14}), (\ref{eqbt12}), we can prove in a similar way as for relation (\ref{eqbt13}) that, for all $\epsilon>0$,
\begin{equation}
\label{eqbt13bis}
\eP \big[\| \ef_i(\ebo)(\eb-\ebo)^t \eM_{2i} (\eb-\ebo) \|_1 \geq \frac{ 6d}{\epsilon} \max   _{1\leq j \leq d} \eE[|\frac{\partial f_i(\ebo)}{\partial \beta_j}|]\big] \leq \frac{\epsilon}{6}.
\end{equation}
\hspace*{0.5cm} For the last term on the right-hand side of (\ref{eq99}), using assumption (A3), we have that, for all $\eb \in {\cal V}_ \eta(\ebo)$, for all $\epsilon>0$, there exists $\epsilon_5>0$, such that $
\eP[\|\eM_{1i}\|_1 \|\eM_{2i} \|_1 \geq \epsilon_5] \leq \epsilon/6$.
 Using this relation, we show similarly, then, for all $\epsilon >0$, there exists $\epsilon_5>0$, such that, 
\begin{equation}
\label{bor1}
\eP \big[ \| \eM_{1i} (\eb-\ebo) (\eb-\ebo)^t \eM_{2i} (\eb-\ebo) \|_1 \geq \epsilon_5 \big]\leq \frac{\epsilon}{6}.
\end{equation}
\hspace*{0.5cm} Choosing
$$M=\sup \Bigg \{ 
\frac{\sigma d }{\sqrt{\epsilon}}\max   _{1\leq j \leq d} \sqrt{6V_{jj}} , \epsilon_2, \frac{6 d }{\epsilon} \max   _{1\leq j,l \leq d} \big \{\eE \big [|\frac{\partial f_i(\ebo)}{\partial \beta_j}\frac{\partial f_i(\ebo)}{\partial \beta_l  }|\big  ] , \eE\big  [|\frac{\partial f_i(\ebo)}{\partial \beta_j}|\big ] \big \}, \epsilon_5 \Bigg \},$$ 
and combining (\ref{eqbt3}), (\ref{eqbt15}), (\ref{eqbt7}), (\ref{eqbt13}), (\ref{eqbt13bis}), (\ref{bor1}) together, lemma yields.
\hspace*{\fill}$\blacksquare$  \\

%%%%%%%%%%%%%%%%%%%%%%%%%%%%%%%%%%%%%%%%%%%%%%%%%%%%%%%%%%%%%%%%%%%%%%%%%%%%
%%%%%%%%%%%%%%%%%%%%%%%%%%%%%%%%%%%%%%%%%%%%%%%%%%%%%%%%%%%%%%%%%%%%%%%%%%%%

\begin{lemma}
\label{lem4}
Under the same assumptions of Theorem \ref{theo2}, we have 
\begin{equation*}
\label{eq12}
\frac{1}{n \theta_{nk}} \sum  _{i \in I} \eg_i(\eb)=O_{\eP}((n \theta_{nk})^{-1/2})+\eV_{1n}^0 (\eb-\ebo)+o_{\eP}(\eb-\eb^0).
\end{equation*}
\end{lemma}
\noindent {\it Proof of Lemma \ref{lem4}}. By the Taylor's expansion up to the order 3 of $\eg_i(\eb)$ at $\eb=\ebo$, we obtain
\begin{eqnarray}
\label{eqd1}
\frac{1}{n\theta_{nk}} \sum   _{i\in I}\eg_i(\eb)&=& \frac{1}{n\theta_{nk}} \sum   _{i\in I}\ef_i(\ebo) \varepsilon_i
 +\frac{1}{2n\theta_{nk}} \sum   _{i\in I}\eff_i(\ebo) (\eb-\ebo) \varepsilon_i \nonumber 
 \\ && 
 -\frac{1}{2n\theta_{nk}} \sum   _{i\in I}\ef_i(\ebo)\ef_i^t (\ebo) (\eb-\ebo)\nonumber 
 \\ && 
  -\frac{1}{6n\theta_{nk}} \sum   _{i\in I}\ef_i(\ebo)(\eb-\ebo)^t \eM_{2i}  (\eb-\ebo) \nonumber 
 \\ &&  
- \frac{1}{4n\theta_{nk}} \sum   _{i\in I} \eff_i(\ebo) (\eb-\ebo)\ef_i^t (\ebo) (\eb-\ebo)
 \\ && 
 -\frac{1}{12n\theta_{nk}} \sum   _{i\in I} \eff_i(\ebo) (\eb-\ebo)(\eb-\ebo)^t  \eM_{2i}  (\eb-\ebo)\nonumber 
 \\ && 
 + \frac{1}{6n\theta_{nk}} \sum   _{i\in I} \eM_{i}  \varepsilon_i - \frac{1}{12n\theta_{nk}} \sum   _{i\in I} \eM_{i}  (\eb-\ebo)^t  \eM_{2i}  (\eb-\ebo), \nonumber
\end{eqnarray} 
 with $\eM_{2i}$ given by Lemma \ref{lem3} and $\eM_{i} = \Bigg( \sum  _{l=1}^d \sum  _{k=1}^d \frac{\partial^2 \ef_i( \eb^{(3)}_{i,kl})}{\partial \beta_k \partial \beta_l} (\beta_k-\beta^0_k)(\beta_l-\beta^0_l)\Bigg)_{1\leq k,l \leq d} $ is  a vector of dimension $(d \times 1)$, where $ \eb^{(3)}_{i,kl} =\ebo+ w_{i,kl}( \eb-\ebo)$, with $ w_{i,kl} \in [0,1]$. \ \\
 For the first term of the right-hand side of (\ref{eqd1}), by the central limit theorem, and the fact that $\eE[\eg_i(\ebo)]=0$, we have 
\begin{equation}
\label{eqd101}
(n \theta_{nk})^{-1} \sum  _{i \in I} \eg_i(\ebo)=O_{\eP}((n \theta_{nk})^{-1/2}).
\end{equation}
For the second term of the right-hand side  of (\ref{eqd1}), by the law of large numbers, the term $ (n\theta_{nk})^{-1} \sum   _{i\in I}\eff_i(\ebo) (\eb-\ebo) \varepsilon_i$ converges almost surely to the expected of $\eff_i(\ebo) (\eb-\ebo) \varepsilon_i$ as $n \rightarrow \infty $. Furthermore, since $\varepsilon_i$ is independent of $\XX_i$ and $\eE[\varepsilon_i]=0$, we have
\begin{equation}
\label{eqd112}
 \frac{1}{n\theta_{nk}} \sum   _{i\in I}\eff_i(\ebo) (\eb-\ebo) \varepsilon_i=o_{\eP}(\eb-\ebo).
\end{equation}
 For the third term of the right-hand side of (\ref{eqd1}), by the law of large numbers and assumption (A4), the term $(n\theta_{nk})^{-1} \sum   _{i\in I}\ef_i(\ebo)\ef_i^t (\ebo) (\eb-\ebo)$ converges almost surely to the expected value  of $\ef_i(\ebo)\ef_i^t (\ebo) (\eb-\ebo)$ as $n \rightarrow \infty $. On the other hand, since $ (n\theta_{nk})^{-1} \sum  _{i \in I} \eff_i(\ebo) \varepsilon_i \overset{a.s}{\longrightarrow} 0$, we have
\begin{equation}
\label{eqd12}
\frac{1}{n\theta_{nk}} \sum   _{i\in I}\ef_i(\ebo)\ef_i^t (\ebo) (\eb-\ebo)=-\eV_{1n}^0 (\eb-\ebo)(1+o_{\eP}(1)).
\end{equation} 
 For the fourth term of the right-hand side of (\ref{eqd1}), by the law of large numbers, using assumption (A3) and the relation (\ref{eqbt12}), we can write $  (6n\theta_{nk})^{-1} \| \sum   _{i\in I}\ef_i(\ebo)$
$(\eb-\ebo)^t \eM_{2i} (\eb-\ebo) \|_1  = O_{\eP}( \| \eb-\ebo \|^2_2)$, which implies 
\begin{equation}
\label{eqd13}
\frac{1}{6n\theta_{nk}} \sum   _{i\in I}\ef_i(\ebo)(\eb-\ebo)^t \eM_{2i} (\eb-\ebo)=o_{\eP}(\eb-\ebo).
\end{equation}
 In the same way, using assumption (A3) and relation (\ref{eqbt12}), we obtain, for the fifth term on the right-hand side of (\ref{eqd1}), that
\begin{equation}
\label{eqd14}
\frac{1}{4n\theta_{nk}} \sum   _{i\in I} \eff_i(\ebo) (\eb-\ebo)\ef_i^t (\ebo) (\eb-\ebo)=o_{\eP}(\eb-\ebo).
\end{equation}
For the sixth term of the right-hand side of (\ref{eqd1}), using the assumption (A3), we have
\begin{equation}
\label{eqd15}
\frac{1}{12n\theta_{nk}} \sum   _{i\in I}\eff_i(\ebo) (\eb-\ebo)(\eb-\ebo)^t \eM_{2i}(\eb-\ebo)=o_{\eP}(\eb-\ebo).
\end{equation}
For $1\leq j \leq d$, and for any fixed $i$, such that $1 \leq i \leq n \theta_{nk}$ , denote by $M_{ij}$ the following random variable designates the j-th component of the vector $\textbf{M}_i$, such that
$$M_{ij}= \sum  _{l=1}^d \sum  _{k=1}^d \frac{\partial^3 f_i(\eb^{(3)}_{i,kl})}{\partial \beta_k \partial \beta_l \partial \beta_j} (\beta_k-\beta^0_k)(\beta_l-\beta^0_l).$$
using assumption (A3), we have with a probability one, $| M_{ij} |\leq C_3 \| \eb-\ebo\|_2^2$.
 Applying Lemma \ref{lem0} (i), we obtain
\begin{equation}
\label{eqMatr}
\|\eM_{i} \|_1 \leq C_3 \| \eb-\ebo\|_2^2.
\end{equation}
 For the term $(6n\theta_{nk})^ {-1} \sum   _{i\in I} \textbf{M}_i \varepsilon_i$, using relations (\ref{epselon}) and (\ref{eqMatr}), we have \ \\
\noindent $(6n\theta_{nk})^{-1} \|\sum   _{i\in I} \textbf{M}_{i} \varepsilon_i\|_1 
\leq (6n\theta_{nk})^{-1} \sum   _{i\in I} \| \textbf{M}_{i}\|_1 |\varepsilon_i| 
\leq C_4(6n\theta_{nk})^{-1} n\theta_{nk} \| \eb-\ebo\|_2^2
= C_4 \| \eb-\ebo\|_2^2$. 
 Then, 
\begin{equation}
\label{Meps}
\frac{1}{6n\theta_{nk}} \sum   _{i\in I} \textbf{M}_{i}\varepsilon_i=o_{\eP}(\eb-\ebo).
\end{equation}
 Finally, for the last term of the right-hand side of (\ref{eqd1}), using assumption (A3) and relation (\ref{eqMatr}), we obtain with probability 1, $ (12n\theta_{nk})^{-1} \| \sum   _{i\in I} \textbf{M}_{i} (\eb-\ebo)^t \eM_{2i} (\eb-\ebo) \|_1
\leq C_5 \| \eb-\ebo\|_2^2$, which gives,
\begin{equation}
\label{eq2013}
 \frac{1}{12n\theta_{nk}} \sum   _{i\in I} \textbf{M}_{i}(\eb-\ebo)^t \eM_{2i}(\eb-\ebo) =o_{\eP}(\eb-\ebo).
\end{equation}
Then, combining relations (\ref{eqd101}), (\ref{eqd112}), (\ref{eqd12}), (\ref{eqd13}), (\ref{eqd14}), (\ref{eqd15}), (\ref{Meps}) and (\ref{eq2013}), we obtain lemma.\hspace*{\fill} \hspace*{\fill}$\blacksquare$ \\

%%%%%%%%%%%%%%%%%%%%%%%%%%%%%%%%%%%%%%%%%%%%%%%%%%%%%%%%%%%%%%%%%%%%%%%%%%%%
%%%%%%%%%%%%%%%%%%%%%%%%%%%%%%%%%%%%%%%%%%%%%%%%%%%%%%%%%%%%%%%%%%%%%%%%%%%%
%%%%%%%%%%%%%%%%%%%%%%%%%%%%%%%%%%%%%%%%%%%%%%%%%%%%%%%%%%%%%%%%%%%%%%%%%%%%

\begin{lemma}
\label{lem1}
 \noindent Under the same assumptions as in Theorem \ref{theo3}, for all $ \varrho
> 0 $, there exist two positive constants $B=B(\varrho)$, $T=T(\varrho)$ such that \ \\
 $ \eP [ \max  _ {\frac{T}{n} \leq \theta_{nk} \leq
1-\frac{T}{n}} (n\theta_{nk} / \log \log n\theta_{nk})^{{1}/{2}} \|
\frac{\hat{\el}(\theta_{nk})}{\min
\{ \theta_{nk}, 1-\theta_{nk} \}} \|_2 > B ] \leq \varrho, \\ 
\eP [ \max  _ {\frac{T}{n} \leq \theta_{nk} \leq
1-\frac{T}{n}} (n\theta_{nk} /  \log \log n\theta_{nk})^{{1}/{2}} \| \hat{\eb}(\theta_{nk})- \ebo \|_2> B ] \leq \varrho, \\ 
\eP [ n^{-{1}/{2}} \max  _ {\frac{T}{n} \leq \theta_{nk}
\leq 1-\frac{T}{n}} n\theta_{nk} \| \frac{\hat{\el}(\theta_{nk})}{\min
\{\theta_{nk}, 1-\theta_{nk}\} } \|_2 > B ] \leq \varrho, \\ 
\eP [ n^{-{1}/{2}} \max  _ {\frac{T}{n} \leq \theta_{nk}
\leq 1-\frac{T}{n}} n\theta_{nk} \| \hat{\eb}(\theta_{nk}) - \ebo \|_2> B ] \leq \varrho.$
\end{lemma}
\noindent {\it Proof of Lemma \ref{lem1}}.  \noindent The proof of this lemma is similar to that of Lemma 1.2.2
of \cite{Csorgo:Horvath:97}. \hspace*{\fill} \hspace*{\fill}$\blacksquare$ \\ \\
In order, to prove Lemma {\ref{lem2}}, we consider 
$$R_k= n \sigma ^{-2}
\theta_{nk}(1-\theta_{nk})(\textbf{W}_{1n}^0-\textbf{W}_{2n}^0)^t\eV^{-1}
(\textbf{W}_{1n}^0-\textbf{W}_{2n}^0).$$  \\
Recall that $\eV \equiv \eE[ \ef(\XX_i,\ebo) \ef^t(\XX_i,\ebo)]$, for all $i=1,...,n$. \\
The results of Lemma {\ref{lem2}} are similar to that of Theorem 1.1.1 of \cite{Csorgo:Horvath:97}. 
\begin{lemma} 
\label{lem2}
Suppose that the assumptions (A1)-(A4) hold. Under the null hypothesis $H_0$, for all $
0 \leq \alpha < 1/2$ we have \ \\
$(i)\,\,\,n^{\alpha} \max  _{\theta_{nk} \in
\Theta_{nk}}[\theta_{nk}(1-\theta_{nk})]^{\alpha} | Z_{nk}(\theta_{nk},
\hat{\el}(\theta_{nk}), \hat{\eb}(\theta_{nk}))-R_k|= O_{\eP}(1)$. \ \\
$(ii)\,\,\max  _{\theta_{nk} \in
\Theta_{nk}}[\theta_{nk}(1-\theta_{nk})]| Z_{nk}(\theta_{nk},
\hat{\el}(\theta_{nk}), \hat{\eb}(\theta_{nk}))-R_k| =O_{\eP}(n^{-1/2}(\log \log n)^{3/2})$.
\end{lemma}
\noindent {\it Proof of Lemma \ref{lem2}}. For the score function $\ephi_{1n}$ of relation $(\ref{eq6})$, the two terms of the right-hand side are replaced by their decomposition obtained by the relations (\ref{eq15}) and (\ref{eq18}). On the other hand, we have $\ephi_{1n}(\theta_{nk},\hat{\el}(\theta_{nk}),\hat{\eb}(\theta_{nk}))= \textbf{0}_d $. Then, we can write  
 $[ \frac{1}{n
\theta_{nk}} \sum  _{i \in I}\eg_i(\ebo)+ \eV_{1n}^0(\hat{\eb}(\theta_{nk})-\ebo)-
\frac {1}{n \theta_{nk}^2}   \sum _{i \in
I}\eg_i(\ebo)\eg^t_i(\ebo)$ $\cdot \hat{\el}(\theta_{nk}) ](1+o_{\eP}(1))
-\eV_{1n}^0 (\eV_{2n}^0)^{-1}[\frac
{1}{n (1-\theta_{nk})^2}\eV_{1n} (\eV_{2n}^0)^{-1}  \sum  
_{j \in J}\eg_j(\ebo)\eg^t_j(\ebo) \hat{\el}(\theta_{nk})+ 
\frac{1}{n(1-\theta_{nk})}$ $ \cdot  \sum   _{j \in
J}\eg_j(\ebo) +\eV_{2n}^0 (\hat{\eb}(\theta_{nk})-\ebo)  ](1+o_{\eP}(1))
=\textbf{0}_d$.  \\
Hence,
\begin{eqnarray*}
\hat{\el}(\theta_{nk}) &=&
\Big (
 \frac{1}{\theta_{nk}} \eD_{1n}^0 + 
\frac{1}{1-\theta_{nk}}(\eV_{1n}^0 (\eV_{2n}^0)^{-1})(\eV_{1n}^0 (\eV_{2n}^0)^{-1})^t
\textbf{D}_{2n}^0 \Big ) ^{-1} \\ &&
\cdot \Big( 
 \frac{1}{n\theta_{nk}}
   \sum   _{i \in I}\eg_i(\ebo)
- \frac{\eV_{1n}^0(\eV_{2n}^0)^{-1}}{n(1-\theta_{nk})}
  \sum   _{j \in
J}\eg_j(\ebo)\Big)(1+o_{\eP}(1)) \\ &&+ o_{\eP}(\hat\eb(\theta_{nk})-\ebo),
\end{eqnarray*}
 with the matrices $\eD^0_{1n}$ and $\eD^0_{2n}$ given by relation (\ref{eqD1D2}).  \ \\
On the other hand, by the law of large numbers, we have $-\eV_{1n}^0 \overset{a.s}{\longrightarrow}
\eV$ and  $-\eV_{2n}^0 \overset{a.s}{\longrightarrow} \eV$. Then, $\eV_{1n}^0 (\eV_{2n}^0
)^{-1}\overset{a.s}{\longrightarrow} I_d$. Always, by the law of large numbers, $\eD_{1n}^0 $ and $\eD_{2n}^0$ converge almost surely to $\sigma^2 \eV$ as $n\rightarrow \infty$. \\
 By Theorem $\ref{theo2}$, we proved that $\hat{\el}(\theta_{nk})=\theta_{nk}
O_p((n\theta_{nk})^{-1/2})$. Then, we obtain 
\begin{equation}
\label{eq49}
 \hat{\el}(\theta_{nk})= \sigma^{-2} \theta_{nk} (1-\theta_{nk})
\eV^{-1} (\textbf{W}_{1n}^0-\textbf{W}_{2n}^0) (1+o_{\eP}(1))+
o_{\eP}(\hat\eb(\theta_{nk})-\ebo).
\end{equation}
The limited development of the statistic
$Z_{nk}(\theta_{nk}, \hat{\el}(\theta_{nk}),\hat{\eb}(\theta_{nk}))$, specified by the relation (\ref{eq5}), in the neighbourhood of $(\el,\eb)=(\textbf{0}_d,\ebo)$ up to order 2, can be written 
\begin{eqnarray}
\label{DLz}
&&\Big [ \frac{2\hat {\el}^t(\theta_{nk})}{\theta_{nk}}  \sum  _{i \in I}\eg_i(\ebo)-  \frac{2\hat{\el}^t(\theta_{nk})}{1-\theta_{nk}} \eV_{1n}^0(\eV_{2n}^0)^{-1} \sum  _{j \in J} \eg_j(\ebo)\Big]   - \Big [ \frac{\hat{\el}^t(\theta_{nk}) }{(1-\theta_{nk})^2} \eV_{1n}^0 (\eV_{2n}^0)^{-1} 
 \nonumber \\ && 
 \cdot \sum  _{j \in J} \eg_j(\ebo)\eg^t_j(\ebo) \eV_{1n}^0 
(\eV_{2n}^0)^{-1}  \hat{\el}(\theta_{nk}) 
+  \frac{\hat{\el}^t(\theta_{nk})}{\theta_{nk}^2}  \sum  _{i \in
I}\eg_i(\ebo)\eg^t_i(\ebo) \hat{\el}(\theta_{nk}) \Big]
 \nonumber \\  && 
 + \Big[ 2\hat{\el}^t(\theta_{nk}) 
  \Big(\frac{1}{\theta_{nk}}
\sum  _{i \in I}\egg_i(\ebo) -\frac{1}{1-\theta_{nk}}  \eV_{1n}^0(\eV_{2n}^0)^{-1} \sum  _{j \in J} \egg_j(\ebo) \Big )  (\hat{\eb}(\theta_{nk})-\ebo) \Big ] 
\nonumber \\  &&  
 - \Big[ 2\hat{\el}^t(\theta_{nk}) \Big(
   \frac{1}{1-\theta_{nk}}  \sum  _{j \in J}\eg_j( \ebo)\frac{\partial (\eV_{1n}(\eb)(\eV_{2n}(\eb))^{-1})}{\partial \eb}\Big ) (\hat{\eb}(\theta_{nk})-\ebo) \Big] 
  \nonumber \\  && 
  + \frac{1}{3 !} \Big [S_1+ 3 S_2 + 3 S_3 +S_4    \Big], 
\end{eqnarray}
%======================================
where \\
$ S_1 = \sum  _{j=1}^d  \sum _{l=1}^d \sum  _{k=1}^d \frac{\partial^3 Z_{nk}(\theta_{nk}, \el^{(1)}_{jkl},  \eb^{(1)}_{jkl})}{\partial \beta_j \partial \beta_k \partial \beta_l } (\hat\beta_j-\beta^0_j) (\hat\beta_k-\beta^0_k)(\hat\beta_l-\beta^0_l)$,  \\
$ S_2 = \sum  _{j=1}^d  \sum  _{l=1}^d \sum  _{k=1}^d \frac{\partial^3 Z_{nk}(\theta_{nk},  \el^{(2)}_{jkl},  \eb^{(2)}_{jkl})}{\partial \lambda_j \partial \lambda_k \partial\beta_l } (\hat\lambda_j) (\hat\lambda_k)(\hat\beta_l-\beta^0_l)$,\\
$ S_3 =  \sum _{j=1}^d  \sum  _{l=1}^d \sum  _{k=1}^d \frac{\partial^3 Z_{nk}(\theta_{nk},  \el^{(3)}_{jkl},  \eb^{(3)}_{jkl})}{\partial \lambda_j \partial \beta_k \partial \beta_l } (\hat\lambda_j) (\hat\beta_k-\beta^0_k)(\hat\beta_l-\beta^0_l)$,\\
$ S_4 = \sum  _{j=1}^d  \sum  _{l=1}^d \sum  _{k=1}^d \frac{\partial^3 Z_{nk}(\theta_{nk}, \el^{(4)}_{jkl},  \eb^{(4)}_{jkl})}{\partial \lambda_j \partial \lambda_k \partial \lambda_l } (\hat\lambda_j) (\hat\lambda_k)(\hat\lambda_l)$,  \\ \\
where, for $1 \leq j \leq d$, $\hat \beta_j$ is the j-th component of $\hat \eb (\theta_{nk})$, and $\hat \lambda_j$ is the \\ j-th component of $\hat \el (\theta_{nk})$. In the expression of $S_1$, $S_2$,$S_3$, $S_4$ we have also, for all $1 \leq j,k,l \leq d$, $\el^{(a)}_{jkl}= u^{(a)}_{jkl} (\hat{\eb}(\theta_{nk})-\ebo)$, and $\eb^{(a)}_{jkl}= \ebo+ v^{(a)}_{jkl} (\hat{\eb}(\theta_{nk})-\ebo)$, with $u^{(a)}_{jkl},v^{(a)}_{jkl} \in [0,1]$ and $a \in \{1,2,3,4\}$. \\
We note that, the derivative $\partial ( \eV_{1n}(\eb)(\eV_{2n}(\eb))^{-1}) / \partial \eb $  is considered term by term. \\
Now, we replace $\hat{\el}(\theta_{nk})$ in the relation (\ref{DLz}) by the value obtained in $(\ref{eq49})$. For the first term of (\ref{DLz}), using notations given by relation (\ref{w1w2}), and the fact that $\eV_{1n}^0 (\eV_{2n}^0
)^{-1}\overset{a.s}{\longrightarrow} I_d$, as $n\rightarrow \infty$, we find that this term is equal to $2n \sigma^{-2}\theta_{nk}(1-\theta_{nk})
(\textbf{W}_{1n}^0-\textbf{W}_{2n}^0)^t\eV^{-1}(\textbf{W}_{1n}^0-\textbf{W}_{2n}^0)+o_{\eP}(\|\hat\eb(\theta_{nk})-\ebo\|_2)$. \\
Similarly, for the second term of (\ref{DLz}), using notations given by (\ref{eqD1D2}), and the fact that $\eD_{1n}^0 $ and $\eD_{2n}^0$
converge to $\sigma^2 \eV$, as $n\rightarrow \infty$, we obtain that this term is equal to $n \sigma^{-2}\theta_{nk}(1-\theta_{nk})(\textbf{W}_{1n}^0-\textbf{W}_{2n}^0)^t\eV^{-1}(\textbf{W}_{1n}^0-\textbf{W}_{2n}^0)+o_{\eP}(\|\hat\eb(\theta_{nk})-\ebo\|_2)$. \\
For the third term of (\ref{DLz}), we know that, $\eV_{1n}^0=(n\theta_{nk})^{-1}
\sum  _{i \in I}\egg_i(\ebo) $, and $\eV_{2n}^0=(n(1-\theta_{nk}))^{-1}
\sum  _{j \in J}\egg_j(\ebo) $. On the other hand, by the law of large numbers, we have $\eV_{1n}^0$ and $\eV_{2n}^0$ converge almost surely to $-\eV$ as $n\rightarrow \infty$, and $\eV_{1n}^0 (\eV_{2n}^0)^{-1}\overset{a.s}{\longrightarrow} I_d$, which implies that the third term of (\ref{DLz}) converge almost surely to zero, as $n\rightarrow \infty$.\\
By the central limit theorem, we have that $(n(1-\theta_{nk}))^{-1}\sum  _{j \in J}\eg_j(\ebo)$ = $O_{\eP}((n(1-\theta_{nk}))^{-1/2})$. Then, the fourth term of (\ref{DLz}) is $o_{\eP}(n \sigma^{-2}\theta_{nk}(1-\theta_{nk})(\textbf{W}_{1n}^0-\textbf{W}_{2n}^0)^t\eV^{-1}(\textbf{W}_{1n}^0-\textbf{W}_{2n}^0)$. \\
For the last term of (\ref{DLz}), using assumptions (A2)-(A4) and by an elementary calculations, we prove that this term is $o_{\eP}(\|\hat\eb(\theta_{nk}) -\ebo\|_2)+o_{\eP}(\|\hat \el(\theta_{nk}) \|_2)+o_{\eP}(\|\hat \el(\theta_{nk})\|_2 \|\hat \eb(\theta_{nk})-\ebo\|_2)$. Combining the obtained results, we obtain 
\begin{eqnarray*}
 && Z(\theta_{nk}, \hat{\el}(\theta_{nk}),\hat{\eb}(\theta_{nk}))=
 n
\sigma^{-2} \theta_{nk}(1-\theta_{nk})
(\textbf{W}_{1n}^0-\textbf{W}_{2n}^0)^t \eV^{-1}
(\textbf{W}_{1n}^0-\textbf{W}_{2n}^0) (1+o_{\eP}(1) ) 
\\  &&  +o_{\eP}(\|\hat \eb(\theta_{nk})-\ebo\|_2)+o_{\eP}(\|\hat \el(\theta_{nk})\|_2) +o_{\eP}(\|\hat \el(\theta_{nk}) \|_2 \|\hat \eb(\theta_{nk}) -\ebo\|_2).
\end{eqnarray*} 
This last relation, together with Lemma {\ref{lem1}} imply Lemma \ref{lem2}.\hspace*{\fill}\hspace*{\fill}$\blacksquare$ \\

\bibliographystyle{plain}

\end{document}